\documentclass[a4paper,11pt,BCOR1cm,parskip]{scrartcl}

\usepackage{german}
\usepackage[iso]{umlaute}
\usepackage{theorem}
\usepackage{amsmath}
\usepackage{amssymb}
\usepackage{makeidx}
\usepackage[matrix,arrow,curve]{xy}
\usepackage{epsfig}
\usepackage{boxedminipage}
\usepackage{epic}
\usepackage{longtable}
\usepackage{ifthen}
\usepackage{tabularx}
\usepackage{parskip}

\newcommand{\R}{\mathrm{I\!R}}
\newcommand{\N}{\mathrm{I\!N}}

\newcommand{\Z}{\mathchoice {\hbox{$\sf\textstyle Z\kern-0.4em
Z$}}{\hbox{$\sf\textstyle Z\kern-0.4em Z$}}{\hbox{$\sf\scriptstyle
Z\kern-0.3em Z$}}{\hbox{$\sf\scriptscriptstyle Z\kern-0.2em Z$}}}

\newcommand{\Q}{\mathchoice {\setbox0=\hbox{$\displaystyle\rm
Q$}\hbox{\raise0.15\ht0\hbox to0pt{\kern0.4\wd0\vrule
height0.8\ht0\hss}\box0}}{\setbox0=\hbox{$\textstyle\rm
Q$}\hbox{\raise0.15\ht0\hbox to0pt{\kern0.4\wd0\vrule
height0.8\ht0\hss}\box0}}{\setbox0=\hbox{$\scriptstyle\rm
Q$}\hbox{\raise0.15\ht0\hbox to0pt{\kern0.4\wd0\vrule
height0.7\ht0\hss}\box0}}{\setbox0=\hbox{$\scriptscriptstyle\rm
Q$}\hbox{\raise0.15\ht0\hbox to0pt{\kern0.4\wd0\vrule
height0.7\ht0\hss}\box0}}}

\newcommand{\id}{\mathrm{id}}

\newcommand{\SO}{\mathrm{SO}}
\newcommand{\End}{\mathrm{End}}

\newcommand{\rk}{\mathrm{rk}}

\newcommand{\eps}{\varepsilon}

\newcommand{\vkap}{\varkappa}

\newcommand{\qmq}[1]{\quad\mbox{#1}\quad}
\newcommand{\Menge}[2]{\{\,#1\,|\,#2\,\}}

\renewcommand{\bigoplus}{\mathop{\bigcirc
  \raisebox{-0.22em}{\hskip-0.53em\hbox{\vrule height2.08ex width0.04em}
  \raisebox{ 0.48em}{\hskip-0.75em\hbox{\vrule height0.04em width 0.8em}}
  \hskip- 0.2em}}}

\newcommand {\g}[2]{\langle #1,#2\rangle}

\newcommand{\sign}{\mathop{\mathrm{sign}}\nolimits}

\newcommand{\Ug}{\mathrm{U}}
\newcommand{\SU}{\mathrm{SU}}

\newcommand{\RE}{\mathop{\mathrm{Re}}\nolimits}
\newcommand{\IM}{\mathop{\mathrm{Im}}\nolimits}

\newcommand{\Eig}{\mathop{\mathrm{Eig}}\nolimits}

\newcommand{\ad}{\mathop{\mathrm{ad}}\nolimits}

\newcommand{\spn}{\mathop{\mathrm{span}}\nolimits}
\newcommand{\tr}{\mathop{\mathrm{tr}}\nolimits}

\newcommand{\liea}{\mathfrak{a}}

\newcommand{\lieg}{\mathfrak{g}}

\newcommand{\liek}{\mathfrak{k}}

\newcommand{\liem}{\mathfrak{m}}

\newcommand{\liet}{\mathfrak{t}}

\newcommand{\liesu}{\mathfrak{su}}

\newcommand{\RP}{\ensuremath{\R\mathrm{P}}}

\newcommand{\beweis}{\begingroup\footnotesize \emph{Proof. }}
\newcommand{\beweisende}{\strut\hfill $\Box$\par\medskip\endgroup}
\newcommand{\Mengegr}[2]{\{\,#1\,{\bigr |}\,#2\,\}}
\newcommand{\wt}{\widetilde}

\newcommand{\C}{\mathchoice {\setbox0=\hbox{$\displaystyle\rm
C$}\hbox{\hbox to0pt{\kern0.4\wd0\vrule
height0.95\ht0\hss}\box0}}{\setbox0=\hbox{$\textstyle\rm C$}\hbox{\hbox
to0pt{\kern0.4\wd0\vrule
height0.95\ht0\hss}\box0}}{\setbox0=\hbox{$\scriptstyle\rm C$}\hbox{\hbox
to0pt{\kern0.4\wd0\vrule
height0.95\ht0\hss}\box0}}{\setbox0=\hbox{$\scriptscriptstyle\rm
C$}\hbox{\hbox to0pt{\kern0.4\wd0\vrule height0.95\ht0\hss}\box0}}}

\theoremstyle{plain} % margin, change

\theorembodyfont{\rmfamily\mdseries\itshape}
\newtheorem{Def}{Definition}[section]
\newtheorem{Prop}[Def]{Proposition}

\theorembodyfont{\rmfamily\mdseries\upshape}

\newtheorem{Remark}[Def]{Remark}
\newtheorem{Remarks}[Def]{Remarks}

\begin{document}
\selectlanguage{english}

\title{Reconstructing the geometric structure \\ of a Riemannian symmetric space \\ from its Satake diagram}
\author{Sebastian Klein${}^1$}
\date{January 24, 2008}
\maketitle
\footnotetext[1]{This work was supported by a fellowship within the Postdoc-Programme of the German Academic Exchange Service (DAAD).}

\abstract{\textbf{Abstract.}} 
The local geometry of a Riemannian symmetric space is described completely by the Riemannian metric and the Riemannian curvature tensor of the space.
In the present article I describe how to compute these tensors for any Riemannian symmetric space from its Satake diagram,
in a way that is suited for the use with computer algebra systems; an example implementation for \textsf{Maple} Version 10 can be found on
\texttt{http://satake.sourceforge.net}. As an example application, the totally geodesic submanifolds of the 
Riemannian symmetric space \,$\SU(3)/\SO(3)$\, are classified.

\bigskip

\textbf{Author's address.} \\
Sebastian Klein \\ Department of Mathematics \\ University College Cork \\ Cork \\ Ireland \\
\texttt{mail@sebastian-klein.de}

\bigskip


\bigskip

\textbf{MS classification numbers:} 53C35 (Primary); 53B20, 17B20, 17-08
% 53C35 Symmetric spaces
% 53B20 Local Riemannian geometry
% 17B20 Simple, semisimple, reductive (super)algebras (roots)
% 17-08 Computational methods

\section{Introduction}
\label{Se:intro}

It is well-known that the behavior of any Riemannian manifold \,$M$\, is influenced strongly by its Riemannian curvature tensor. To give
just two examples, the ``spreading'' of the geodesics,
as measured by the Jacobi fields, and the existence of a totally geodesic submanifold tangential to a given subspace of a tangent space
are expressed in terms of the curvature tensor field on the manifold \,$M$\,.

Especially for (locally) Riemannian symmetric spaces, the control exerted on the local geometry of the manifold
by the Riemannian curvature, together with the Riemannian metric,
is total: If \,$M$\, and \,$N$\, are two such spaces, and there exists a linear isometry \,$T_pM \to T_qN$\, which transports the curvature tensor of \,$M$\, at \,$p$\,
into the curvature tensor of \,$N$\, at \,$q$\,, then \,$M$\, and \,$N$\, are already locally isometric to each other; note that it suffices to consider the curvature
in a single point of \,$M$\, and \,$N$\, because the curvature tensor field is parallel in this situation.
This shows that for a Riemannian symmetric space \,$M$\,,
the local geometry of \,$M$\, is described completely by two tensors on a single tangent space \,$T_pM$\,: the inner product given by the Riemannian
metric at \,$p$\,, and the curvature tensor at \,$p$\,. 
Viewed in this way, the study of the local geometry of a Riemannian symmetric space reduces to a purely algebraic problem, namely to the study of
these two tensors on the tangent space \,$T_pM$\,. Thus we will call these two tensors the ``fundamental geometric tensors'' of \,$M$\,.

One very important example for the control of the geometry of a Riemannian symmetric space \,$M$\, by its curvature tensor \,$R$\, is the following result,
which permits the classification of the totally geodesic submanifolds of \,$M$\,: A linear subspace \,$U \subset T_pM$\, is the tangent space of a totally
geodesic submanifold of \,$M$\, if and only if \,$U$\, is a Lie triple system, i.e.~if \,$R(u,v)w\in U$\, holds for all \,$u,v,w \in U$\,. 

In view of the above, it is very desirable to have representations of the fundamental geometric tensors,
especially the curvature tensor, available for study for every Riemannian symmetric space
\,$M=G/K$\,. The well-known formula \,$R(u,v)w = -[[u,v],w]$\, relating the curvature tensor \,$R$\, of \,$M$\, to the Lie bracket of the Lie algebra \,$\lieg$\,
of the transvection group \,$G$\, of \,$M$\, lets one calculate \,$R$\, relatively easily if \,$G$\, is a classical group (then \,$\lieg$\, is a matrix Lie algebra,
with the Lie bracket being simply the commutator of matrices), but not so easily if \,$G$\, is one of the exceptional Lie groups, because then the explicit
description of \,$\lieg$\, as a matrix algebra is too unwieldy to be useful generally.

Therefore, in the present paper, I will describe another representation of the fundamental geometric tensors of any Riemannian symmetric space of compact type,
based on the root space decomposition of the Lie algebra of its transvection group.
This representation is especially suited for the use with a computer algebra system.
I have implemented the algorithms and equations given here as a \textsf{Maple} package, which can be found at \texttt{http://satake.sourceforge.net}\,.
In a forthcoming paper, I will use this presentation to classify the totally geodesic submanifolds in the exceptional Riemannian symmetric spaces of rank~\,$2$\,
(based on similar methods as my classification in the \,$2$-Grassmannians, see \cite{Klein:2007-claQ} and \cite{Klein:2007-claG2}). 

As information about the symmetric space concerned, we will require only the Satake diagram of that space, i.e.~the Dynkin diagram
of the Lie algebra of the transvection group, ``annotated'' with further information describing the symmetric structure of the space, see for example 
\cite{Loos:1969-2}, Section~VII.3.3, p.~132ff.. The Satake diagrams are well-known and tabulated in the literature (for example, in \cite{Loos:1969-2}, p.~147f.)
for every irreducible Riemannian symmetric space. It is a well-known fact that the Satake diagram already determines the local structure
of the Riemannian symmetric space; however, it turns out that for the actual reconstruction of the fundamental geometric tensors in a sufficiently explicit way,
some new work needs to be done.

Our consideration is based on the following well-known construction: 
Let us consider a Riemannian symmetric space \,$M=G/K$\, of compact type. 
Then the symmetric structure of \,$M$\,
induces an involutive automorphism \,$\sigma$\, on the Lie algebra \,$\lieg$\, of the transvection group \,$G$\, of \,$M$\,, and \,$\sigma$\,
gives rise to the decomposition \,$\lieg = \liek \oplus \liem$\,, where \,$\liek = \Eig(\sigma,1)$\, is the Lie algebra of the isotropy group \,$K$\,
and \,$\liem = \Eig(\sigma,-1)$\, is a linear subspace of \,$\lieg$\, which is canonically isomorphic to the tangent space \,$T_oM$\, at the ``origin point''
\,$o := eK \in G/K=M$\,, and if we identify \,$T_oM$\, with \,$\liem$\, in that way, then on each irreducible factor of \,$M$\,, the Riemannian metric of \,$M$\,
is a constant multiple of the Killing form of \,$\lieg$\,, and the curvature tensor \,$R$\, of \,$M$\, at \,$o$\, is given by
the formula \,$R(u,v)w = -[[u,v],w]$\, for \,$u,v,w \in \liem$\,. 

To reconstruct the curvature tensor and the inner product on \,$T_oM$\, from the Satake diagram, we proceed in the following way:
It is well-known that the Dynkin diagram of \,$\lieg$\, (which can be read off the Satake diagram of \,$M$\,) determines uniquely the
root system \,$\Delta$\, of \,$\lieg$\,; an algorithm for the reconstruction of \,$\Delta$\, form the Dynkin diagram is given in Section~\ref{Se:root}. 
Because the root spaces \,$\lieg^{\C}_\alpha$\, of the complexification \,$\lieg^{\C}$\, of \,$\lieg$\,
are complex-1-dimensional, the action of the Lie bracket on each \,$\lieg^{\C}_\alpha \times \lieg^{\C}_\beta$\, is already determined 
up to a constant by the structure of \,$\Delta$\, via the 
relation \,$[\lieg^{\C}_\alpha,\lieg^{\C}_\beta]=\lieg^{\C}_{\alpha+\beta}$\,; more specifically, if we choose a non-zero
\,$X_\alpha \in \lieg^{\C}_\alpha$\, for each \,$\alpha\in\Delta$\,, then there exist constants \,$c_{\alpha,\beta}\in\C$\, so that 
\,$[X_\alpha,X_\beta] = c_{\alpha,\beta}\,X_{\alpha+\beta}$\, holds for every \,$\alpha,\beta\in\Delta$\, with \,$\alpha+\beta\in\Delta$\,. 
As has been discovered by \textsc{Weyl} and \textsc{Chevalley} (see \cite{Weyl:1933} and \cite{Chevalley:1955}), 
there exists a way to choose the vectors \,$X_\alpha$\, in such a way that
the constants \,$c_{\alpha,\beta}$\, become real, and can \emph{up to sign} be computed by a simple formula dependent only on the structure
of the root system \,$\Delta$\, (see also Proposition~\ref{P:c:c}(g) of the present paper).
We call a system \,$(X_\alpha)_{\alpha\in\Delta}$\, chosen in accordance with this a \emph{Chevalley basis} of \,$\lieg^{\C}$\,.
In Section~\ref{Se:c} we show how a Chevalley basis can further be adapted to the position of the compact Lie algebra \,$\lieg$\, within \,$\lieg^{\C}$\,, 
and then give an algorithm to compute the constants \,$c_{\alpha,\beta}$\, corresponding to such a further adapted basis, 
including their sign, thereby recovering the Lie algebra structure of \,$\lieg$\, completely.

Section~\ref{Se:sigma} then discusses the action of the involutive automorphism \,$\sigma:\lieg\to\lieg$\,, which describes the symmetric structure of \,$M$\,,
on \,$\lieg$\,. The action of \,$\sigma$\, on the Cartan subalgebra (spanned by the root vectors) is already well-known from the works of \textsc{Satake}
(indeed, its description was what induced Satake to introduce what is now known as the Satake diagram), so it remains for Section~\ref{Se:sigma}
to describe how \,$\sigma$\, acts on a Chevalley basis of \,$\lieg^{\C}$\,, and thus on
the root spaces of \,$\lieg$\,. By knowing \,$\sigma$\,, we then know the splitting \,$\lieg = \liek \oplus \liem$\,,
and by also knowing the Lie bracket structure of \,$\lieg$\,, we are then able to calculate the curvature tensor via the formula \,$R(u,v)w = -[[u,v],w]$\,
for \,$u,v,w \in \liem$\,. Because we moreover know how the Killing form evaluates for members of the Chevalley basis we chose, 
we can also express the Killing form on \,$\liem$\,,
and hence the Riemannian metric of \,$M$\, (which is a multiple of the Killing form on each irreducible factor of \,$M$\,) acting on \,$\liem$\,. The resulting formulas, which
describe the fundamental geometric tensors of \,$M$\,, are given in Section~\ref{Se:geom}.

Finally, to illustrate the usefulness of the presentation of the fundamental geometric tensors given herein, I apply it in Section~\ref{Se:AI(2)} to 
classify the Lie triple systems, and thus the totally geodesic submanifolds, of a specific Riemannian symmetric space. Although the presentation
was developed mainly with the exceptional Riemannian symmetric spaces in mind, as explained above, for the sake of simplicity I here investigate
the classical Riemannian symmetric space \,$\SU(3)/\SO(3)$\,. As mentioned above,
I will use the same methods in a forthcoming paper to classify the totally geodesic submanifolds in all exceptional Riemannian symmetric spaces of rank \,$2$\,. 

\bigskip

I have produced an example implementation of the presentation of the fundamental geometric tensors described in this article as a package for
\textsf{Maple} Version 10. This implementation can be downloaded from \verb!http://satake.sourceforge.net!. The worksheet and the corresponding
technical documentation also accompanies the version of the present paper posted on \verb!http://www.arxiv.org!.

%\bigskip
\newpage

The results of the present paper were obtained by me while working at the University College Cork under the advisorship of Professor J.~Berndt.
I would like to thank him for his dedicated support and guidance, as well as his generous hospitality.

\section{Reconstructing the root system}
\label{Se:root}

As described in the Introduction, the first step in the reconstruction of the fundamental geometric tensors of a Riemannian symmetric space (of compact type)
is the reconstruction of the root system of the Lie algebra of its transvection group from its Dynkin diagram.

For this purpose, we 
let \,$\lieg$\, be a compact real Lie algebra, i.e.~the Killing form \,$\vkap: \lieg \times \lieg \to \R,\; (X,Y) \mapsto \tr(\ad(X) \circ \ad(Y))$\, 
of \,$\lieg$\, is negative definite. We fix a Cartan subalgebra
(i.e.~a maximal abelian subalgebra) \,$\liet$\, of \,$\lieg$\,. For any \,$\C$-linear form \,$\alpha\in  (\liet^{\C})^*$\, on the complexification
\,$\liet^{\C} := \liet \otimes_{\R} \C$\, of \,$\liet$\, we consider 
$$ \lieg_\alpha := \Menge{X \in \lieg}{\forall H \in \liet: \ad(H)^2 X = \alpha(H)^2 X} \;; $$
if we have \,$\alpha \neq 0$\, and \,$\lieg_\alpha \neq \{0\}$\,, \,$\alpha$\, is called a \emph{root} of \,$\lieg$\, (with respect to \,$\liet$\,)
and \,$\lieg_\alpha$\, is called the \emph{root space} corresponding to \,$\alpha$\,;
the set \,$\Delta$\, of all roots of \,$\lieg$\, is called the \emph{root system} of \,$\lieg$\,.

Because the Lie algebra \,$\lieg$\, is compact, the roots of \,$\lieg$\, are purely imaginary on \,$\liet$\,, i.e.~each \,$\alpha \in \Delta$\, is of the
form \,$\alpha = i\,\alpha'$\,, where \,$\alpha' \in \liet^*$\, is a real linear form on \,$\liet$\,. 
It is clear that \,$\lieg_{-\alpha} = \lieg_\alpha$\, holds, and therefore we have
\,$-\alpha \in \Delta$\, if and only if \,$\alpha \in \Delta$\,. Thus, if we fix \,$H_0 \in \liea$\, so that \,$\alpha(H_0) \neq 0$\,
holds for all \,$\alpha \in \Delta$\, (such a \,$H_0$\, exists because \,$\Delta$\, is finite), then the subset
$$ \Delta_+ := \Menge{\alpha \in \Delta}{\alpha(H_0) \in i\,\R_+} $$
of \,$\Delta$\, (called the set of \emph{positive roots} with respect to \,$H_0$\,) satisfies 
\,$\Delta_+ \cup (-\Delta_+) = \Delta$\, and \,$\Delta_+ \cap (-\Delta_+) = \varnothing$\,; with respect to it we have the
\emph{root space decomposition} of \,$\lieg$\,
\begin{equation}
\label{eq:root:decomp-g}
\lieg = \liet \;\oplus\; \bigoplus_{\alpha \in \Delta_+} \lieg_\alpha \; .
\end{equation}

A root \,$\alpha \in \Delta_+$\, is called \emph{simple}, if it is not the sum of two positive roots. We denote the set of simple roots in \,$\Delta_+$\,
by \,$\Pi$\,; it is a basis of \,$i\liet^*$\,. If \,$\beta \in \Delta$\, is an arbitrary root, then the coefficients \,$k_\alpha$\, in the unique representation
\,$\beta = \sum_{\alpha\in\Pi}k_\alpha\cdot\alpha$\, are all integers; moreover they are either all \,$\geq 0$\, (this is the case if and only if \,$\beta \in \Delta_+$\,)
or all \,$\leq 0$\, (this is the case if and only if \,$\beta \in -\Delta_+$\,). (See \cite{Knapp:2002}, Proposition~2.49, p.~155.) Therefore the number
$$ \ell(\beta) := \sum_{\alpha\in\Pi} k_\alpha \;, $$
called the \emph{level} of \,$\beta$\,, is an integer \,$\neq 0$\,; we have \,$\ell(\beta)>0$\, if and only if \,$\beta \in \Delta_+$\, holds,
and \,$\ell(\beta)=1$\, if and only if \,$\beta$\, is simple.
Clearly the level is additive, i.e.~we have for any \,$\alpha,\beta \in \Delta$\, with \,$\alpha+\beta \in \Delta$\,
\begin{equation}
\ell(\alpha+\beta) = \ell(\alpha) + \ell(\beta)\; . 
\end{equation}

Next, we use the Killing form \,$\vkap$\, of \,$\lieg$\, to induce an inner product on \,$i\liet$\, resp.~on \,$i\liet^* = \spn_{\R}(\Delta)$\,:
Because the Lie algebra \,$\lieg$\, is compact, 
\,$\vkap$\, is negative definite on \,$\lieg$\,, and therefore its complexification, which we again denote by 
\,$\vkap: \lieg^{\C} \times \lieg^{\C} \to \C$\,, is non-degenerate; its restriction \,$\g{\,\cdot\,}{\,\cdot\,} := \vkap|(i\liet \times i\liet)$\,
is a positive definite inner product on \,$i\liet$\,. It follows that for any \,$\alpha \in (i\liet)^*$\, there exists one and only one
vector \,$\alpha^\sharp \in i\liet$\, so that \,$\alpha = \vkap(\alpha^\sharp,\,\cdot\,)|i\liet$\, holds. By pulling back the inner product
of \,$i\liet$\, with the linear isomorphism \,$(i\liet)^* \to i\liet,\;\alpha\mapsto\alpha^\sharp$\,, we get an inner product on \,$(i\liet)^*$\,,
which we will also denote by \,$\g{\,\cdot\,}{\,\cdot\,}$\,. We will also use the associated norm \,$\|v\| := \sqrt{\g{v}{v}}$\, for \,$v \in i\liet$\,
resp.~for \,$v\in (i\liet)^*$\,. 

It is a well-known fact that for any two roots \,$\alpha,\beta\in\Delta$\,, the quantity
\,$n_{\alpha,\beta} := 2 \tfrac{\g{\alpha}{\beta}}{\|\alpha\|^2}$\, relating the lengths and the angles of the roots
can only attain the discrete values \,$0,\pm 1, \pm 2, \pm 3$\, (see, for example, \cite{Knapp:2002}, Proposition~2.48(c), p.~153);
if \,$\alpha,\beta$\, are simple with \,$\alpha\neq\beta$\,, we moreover have \,$n_{\alpha,\beta}\cdot n_{\beta,\alpha} \in \{0,1,2,3\}$\,.
The matrix \,$N := (n_{\alpha,\beta})_{\alpha,\beta \in \Pi}$\, is called the \emph{Cartan matrix} of \,$\Delta_+$\, (or of \,$\lieg$\,); it is known
that up to conjugation with a permutation matrix, it does not depend on the choices involved in obtaining \,$\Delta_+$\,. Moreover, \,$N$\, is regular,
and if \,$\lieg$\, is simple, \,$N$\, does not have a non-trivial block diagonal form. Finally the \emph{Dynkin diagram} of \,$\lieg$\, is obtained
by drawing a node for each simple root of \,$\lieg$\,; the nodes corresponding to \,$\alpha,\beta\in\Pi$\, are connected by the number of
lines indicated by \,$n_{\alpha,\beta}\cdot n_{\beta,\alpha} \in \{0,1,2,3\}$\,; if \,$n_{\alpha,\beta}\cdot n_{\beta,\alpha} \in \{2,3\}$\, holds, 
then \,$\alpha$\, and \,$\beta$\, are of unequal length, and we indicate the relation by drawing an arrow tip pointing from the longer to the shorter root.

\bigskip

We are now able to describe how the root system \,$\Delta$\, and the Killing form of \,$\lieg$\, (up to a constant factor on each simple ideal of \,$\lieg$\,)
are reconstructed from the Dynkin diagram of \,$\lieg$\,: First, from the Dynkin diagram, we can reconstruct the Cartan matrix
\,$N = (n_{\alpha,\beta})_{\alpha,\beta\in\Pi}$\, (where \,$\Pi$\, denotes the set of simple roots, i.e.~of nodes of the Dynkin diagram): Let \,$\alpha,\beta\in\Pi$\,
be given, then \,$n_{\alpha,\beta}$\, is obtained in the following way: If \,$\alpha=\beta$\, holds, we obviously have \,$n_{\alpha,\beta} = 2$\,.
Otherwise, let \,$k$\, denote the number of lines connecting \,$\alpha$\, and \,$\beta$\, in the Dynkin diagram, and for \,$k\in\{2,3\}$\, suppose that
\,$\alpha,\beta$\, are arranged in such a way that \,$\alpha$\, is the longer of the two roots (as indicated in the Dynkin diagram). Then for \,$k=0$\,
we have \,$n_{\alpha,\beta} = n_{\beta,\alpha} = 0$\,, whereas for \,$k \in \{1,2,3\}$\, we have \,$n_{\alpha,\beta}=1$\, and \,$n_{\beta,\alpha}=k$\,. 

From the Cartan matrix we can reconstruct the relative lengths of the simple roots within each simple ideal of the semisimple Lie algebra \,$\lieg$\,:
If \,$\alpha,\beta\in\Pi$\, are simple roots with \,$\g{\alpha}{\beta}\neq 0$\,, then we have \,$\tfrac{\|\beta\|}{\|\alpha\|}
= \sqrt{\tfrac{n_{\beta,\alpha}}{n_{\alpha,\beta}}}$\,; because any two simple roots \,$\alpha,\beta$\, in the same simple ideal of \,$\lieg$\, are connected by
a chain of simple roots \,$\alpha=\gamma_1,\gamma_2,\dotsc,\gamma_{k-1},\gamma_k=\beta$\, with \,$\g{\gamma_j}{\gamma_{j+1}}\neq 0$\, (the Dynkin diagram
of each simple ideal is connected), we thereby know the relative lengths of any two such simple roots. By arbitrarily fixing the length of one simple
root in each simple ideal of \,$\lieg$\, (this corresponds to the choice of a factor for the metric on each simple ideal of \,$\lieg$\,), we then know
the length of each simple root of \,$\lieg$\,. Thereby we also
know the inner product between any two simple roots \,$\alpha,\beta \in \Pi$\,: \,$\g{\alpha}{\beta} = \tfrac12 \, \|\alpha\|^2 \,n_{\alpha,\beta}$\,.
Because \,$\Pi$\, is a basis of \,$i\liet^*$\,, this relationship permits us to reconstruct the inner product \,$\g{\,\cdot\,}{\,\cdot\,}$\, on 
all of \,$i\liet^*$\,.

We now state the algorithm for the reconstruction of \,$\Delta_+$\,. The algorithm in fact reconstructs the roots ordered by level, i.e.~it constructs
the sets \,$\Delta_j := \Menge{\alpha \in \Delta_+}{\ell(\alpha)=j}$\, for all \,$j \in \{1,\dotsc,L\}$\,, where \,$L$\, is the maximal level occurring in \,$\Delta_+$\,. 

\medskip

\begin{enumerate}
\item[\textbf{(R1)}] [Initialization.]
Let \,$\Delta_1$\, be the set of simple roots. 
Let \,$\Delta_\ell := \varnothing$\, for all \,$\ell\geq 2$\, which are needed below.
\item[\textbf{(R2)}] [Iterate on level.] 
Let \,$\ell := 1$\,. Iterate steps (R3)--(R8) until the condition given in (R8) is satisfied.
\item[\textbf{(R3)}] [Iterate on roots.] 
Iterate steps (R4)--(R7) for all \,$\beta \in \Delta_\ell$\, and \,$\alpha \in \Delta_1$\,.
\item[\textbf{(R4)}] [Skip, if the \,$\alpha$-string through \,$\beta$\, has already been generated, or if \,$\beta=\alpha$\,.] 
If we have \,$\beta=\alpha$\,, go to step (R7). 
If we have \,$\ell \geq 2$\, and \,$\beta-\alpha \in \Delta_{\ell-1}$\,, also go to step (R7).
\item[\textbf{(R5)}] [Determine the length of the \,$\alpha$-string through \,$\beta$\,.]
Put \,$q := -2 \tfrac{\g{\beta}{\alpha}}{\|\alpha\|^2}$\,, where the inner product and the norm is calculated as described above.
\item[\textbf{(R6)}] [Add the \,$\alpha$-string through \,$\beta$\,.]
If \,$q \geq 1$\, holds: For each \,$k \in \{1,\dotsc,q\}$\,, put \,$\beta+k\alpha$\, into \,$\Delta_{\ell+k}$\,, if it is not already present there.
\item[\textbf{(R7)}] (End of the body of the loop started in (R3).)
\item[\textbf{(R8)}] Increase \,$\ell$\, by \,$1$\,. The loop started in (R2) (and the algorithm) ends, if we now have \,$\Delta_\ell = \varnothing$\,. Then \,$\ell-1$\,
is the maximal level \,$L$\, of \,$\Delta$\,. 
\end{enumerate}

\medskip

{\footnotesize
\emph{Proof for the correctness of the algorithm.}
We base the proof of the correctness of the algorithm on the concept of a \emph{string} of roots: For \,$\beta \in \Delta \cup \{0\}$\, and \,$\alpha \in \Delta$\,,
the \emph{\,$\alpha$-string through \,$\beta$\,} is the set of all \,$\beta+k\alpha$\,, \,$k \in \Z$\,, so that \,$\beta+k\alpha \in \Delta \cup \{0\}$\, holds.
The \,$\alpha$-string through \,$\beta$\, does not have any gaps, in other words there exist numbers \,$p,q \in \N_0$\, so that this string
is equal to \,$\Menge{\beta+k\alpha}{-p\leq k \leq q}$\,; moreover we have
\begin{equation}
\label{eq:posrts:stringpq}
p-q = 2 \frac{\g{\beta}{\alpha}}{\|\alpha\|^2} \; ,
\end{equation}
see for example \cite{Knapp:2002}, Proposition~2.48(g), p.~153.

Because a root is of level \,$1$\, if and only if it is simple, \,$\Delta_1$\, is given the correct value in step (R1), and \,$\Delta_1$\, is not modified
thereafter in the course of the algorithm. We shall next show by induction that after the iteration of the loop (R2)--(R8) for a given value of \,$\ell$\,
is completed, \,$\Delta_{\ell+1}$\, is exactly the set of all positive roots of level \,$\ell+1$\,. 

First we show that the elements inserted into \,$\Delta_\ell$\, are indeed roots: 
Note that the only insertions into any \,$\Delta_{j}$\, occur in step (R6), and 
this insertion instruction is reached only with \,$\beta$\, being a root of level \,$\ell$\, and \,$\alpha$\, being a simple root. 
Moreover, if \,$\beta-\alpha$\, is a root, then it is of level \,$\ell-1$\,, and by the induction hypothesis
we know that at this stage \,$\Delta_{\ell-1}$\, already contains all roots of that level. Therefore the condition in step (R4) ensures that 
steps (R5)/(R6) are reached only if \,$\beta-\alpha$\, is not a root.
In that case, the \,$\alpha$-string through \,$\beta$\, therefore is of the
form \,$\Menge{\beta+k\alpha}{0\leq k \leq q}$\, with some \,$q \in \N_0$\,, and Equation~\eqref{eq:posrts:stringpq} shows that this \,$q$\, is the number calculated
in step (R5). Therefore the elements \,$\beta+k\alpha$\, inserted into \,$\Delta_{\ell+k}$\, in step (R6) constitute exactly the \,$\alpha$-string through \,$\beta$\, 
and are thus in particular roots; moreover the level of \,$\beta+k\alpha$\, equals \,$\ell+k$\,, and therefore this root is inserted into the correct
set \,$\Delta_{\ell+k}$\,. 

We now show that indeed all roots of level \,$\ell+1$\, have been inserted into \,$\Delta_{\ell+1}$\, by the end of the iteration \,$\ell$\,. 
For this, let \,$\beta$\, be any root of level \,$\ell+1$\,. 
It is known that there exists a simple root \,$\alpha \in \Delta_1$\, so that \,$\beta-\alpha$\, is also a root.
Consider the \,$\alpha$-string through \,$\beta$\,,
\,$\Menge{\beta+k\alpha}{-p\leq k\leq q}$\,; because \,$\beta-\alpha$\, is a root, we have \,$p \geq 1$\,. The root of minimal level in this string,
\,$\beta-p\alpha$\,, is of level \,$\ell+1-p \leq \ell$\,; by induction it has therefore already been generated as a member of \,$\Delta_{\ell+1-p}$\, at 
the beginning of the loop iteration \,$\ell$\,. As shown in the first part of the proof, in the loop iteration \,$\ell-p$\,, with the generation
of \,$\beta-p\alpha$\,, all members of the \,$\alpha$-string through \,$\beta-p\alpha$\,, have been added to the appropriate \,$\Delta_{(\ell+1-p)+k}$\,;
in particular \,$\beta$\, itself has been added to \,$\Delta_{\ell+1}$\,. 
\strut\hfill $\Box$

}

\section{Determining the Lie bracket}
\label{Se:c}

Our next task is to reconstruct the action of the Lie bracket of \,$\lieg$\, from its root system \,$\Delta$\, (which was obtained from the Dynkin diagram of \,$\lieg$\,
in Section~\ref{Se:root}). To do so, we need to take a detour into the complex setting.

We again let \,$\lieg$\, be a compact Lie algebra, and use the notations of the preceding section.
We consider the complexification  \,$\lieg^{\C}$\, of \,$\lieg$\,;
via the complexification of the Lie bracket of \,$\lieg$\,, \,$\lieg^{\C}$\, becomes a complex semisimple Lie algebra. 
It should be noted that the Killing form of \,$\lieg^{\C}$\, equals the complexification \,$\vkap$\, of the Killing form of \,$\lieg$\,. 
The complexification \,$\liet^{\C}$\, of the Cartan subalgebra \,$\liet$\, of \,$\lieg$\, is a Cartan subalgebra of \,$\lieg^{\C}$\,, and we put
for any \,$\alpha \in (\liet^{\C})^*$\,
\begin{align*}
\lieg_\alpha^{\C} & := \Menge{X \in \lieg^{\C}}{\forall H \in \liet^{\C}: \ad(H)X = \alpha(H)X} \\
& = \Menge{X \in \lieg^{\C}}{\forall H \in \liet: \ad(H)X = \alpha(H)X} \,.
\end{align*}
Then the root system \,$\Menge{\alpha \in (\liet^{\C})^* \setminus \{0\}}{\lieg_\alpha^{\C} \neq \{0\}}$\, of \,$\lieg^{\C}$\, equals the root system \,$\Delta$\,
of \,$\lieg$\,, and we have the root space decomposition
$$ \lieg^{\C} = \liet^{\C} \;\oplus\; \bigoplus_{\alpha \in \Delta} \lieg_\alpha^{\C}\;. $$
It is a well-known fact that all these root spaces \,$\lieg_\alpha^{\C}$\,, \,$\alpha\in\Delta$\,, are complex-1-dimensional.

It was famously described by \textsc{Weyl} and \textsc{Chevalley} (see \cite{Weyl:1933} and \cite{Chevalley:1955}) how to choose bases of the root spaces
\,$\lieg_\alpha^{\C}$\, (they consist of only one vector each, because the root spaces are complex-1-dimensional) which are adapted to the Lie bracket of \,$\lieg^{\C}$\,
in the best possible way.  We here cite the result in the form given in \cite{Knapp:2002}, \S VI.1.

\begin{Def}
\label{D:c:chevalley}
A family of vectors \,$(X_\alpha)_{\alpha \in \Delta}$\, is called a \emph{Chevalley basis} of \,$\lieg^{\C}$\,, if we have \,$X_\alpha \in \lieg^{\C}_\alpha$\, for every
\,$\alpha \in \Delta$\, and if there exists a family of \emph{real} numbers \,$(c_{\alpha,\beta})_{\alpha,\beta \in \Delta}$\,, called the \emph{Chevalley constants}
corresponding to \,$(X_\alpha)$\,, so that for all \,$\alpha,\beta \in \Delta$\, we have
\begin{equation}
\label{eq:c:XaXb}
[X_\alpha,X_\beta] = \begin{cases}
c_{\alpha,\beta}\,X_{\alpha+\beta} & \text{if \,$\alpha+\beta\in\Delta$\,} \\
\alpha^\sharp & \text{if \,$\alpha+\beta=0$\,} \\
0 & \text{otherwise}
\end{cases} \; \;,
\end{equation}
and 
\begin{equation}
\label{eq:c:c-a-b}
c_{-\alpha,-\beta} = -c_{\alpha,\beta} \; .
\end{equation}
For formal reasons we put \,$c_{\alpha,\beta} := 0$\, wherever \,$\alpha,\beta \in \Delta$\, with \,$\alpha+\beta \not\in \Delta$\, and \,$\beta\neq -\alpha$\,. 
\end{Def}

It should be noted that the Chevalley constants do depend on the choice of the Chevalley basis. However, their squares are uniquely determined
by the structure of the Lie algebra (as Proposition~\ref{P:c:c}(g) below shows), therefore the transition from one Chevalley basis to another
can change the corresponding Chevalley constants only in sign. The specific transformation behavior of the Chevalley constants is described in 
Proposition~\ref{P:c:chevtrafo} below.

\begin{Prop}
\,$\lieg^{\C}$\, has a Chevalley basis.
\end{Prop}

\beweis
See \cite{Knapp:2002}, Theorem~6.6, p.~351.
\beweisende

\begin{Prop}
\label{P:c:c} 
Let \,$(X_\alpha)$\, be a Chevalley basis of \,$\lieg^{\C}$\, and \,$(c_{\alpha,\beta})$\, be the corresponding Chevalley constants.
Suppose \,$\alpha,\beta,\gamma,\delta \in \Delta$\,. Then we have
\begin{enumerate}
\item 
\,$c_{\beta,\alpha} = -c_{\alpha,\beta}$\,.
\item 
\,$\vkap(X_\alpha,X_{-\alpha}) = 1$\,, where \,$\vkap$\, is the Killing form of \,$\lieg^{\C}$\,. 
\item
\,$\overline{X_\alpha} = a\cdot X_{-\alpha}$\,, with \,$a:= \vkap(X_\alpha,\overline{X_\alpha}) < 0$\,. 
\item
We have \,$\lieg_\alpha = \Menge{V_\alpha(c)}{c \in \C}$\,, where we put
\,$V_\alpha(c) := \tfrac{1}{\sqrt{2}}\,(c\,X_\alpha + \overline{c}\,\overline{X_{\alpha}})$\, for  \,$c \in \C$\,.
\item
Suppose \,$\alpha+\beta+\gamma=0$\,. Then we have \,$c_{\alpha,\beta} = c_{\beta,\gamma} = c_{\gamma,\alpha}$\,.
\item
Suppose \,$\alpha+\beta+\gamma+\delta=0$\, and that none of the roots \,$\alpha,\beta,\gamma,\delta$\, is the negative of one of the others.
Then we have \,$c_{\alpha,\beta}\,c_{\gamma,\delta} + c_{\beta,\gamma}\,c_{\alpha,\delta} + c_{\gamma,\alpha}\,c_{\beta,\delta} = 0$\,.
\item
We have
$$c_{\alpha,\beta}^2 = \frac{q\cdot(1+p)}{2}\cdot \|\alpha\|^2 \;, $$
where \,$\Menge{\beta+k\alpha}{-p \leq k \leq q}$\, is the \,$\alpha$-string through \,$\beta$\,; note that this implies that we have \,$c_{\alpha,\beta} \neq 0$\,
if \,$\alpha+\beta\in \Delta$\, holds. 
\end{enumerate}
\end{Prop}

\beweis
\emph{(a)} is obvious. 
\emph{For (b),} we have \,$\alpha(\alpha^\sharp)\cdot \vkap(X_\alpha,X_{-\alpha}) = \vkap(\ad(\alpha^\sharp)X_\alpha,X_{-\alpha}) = \vkap(\alpha^\sharp,[X_\alpha,X_{-\alpha}])
\overset{\eqref{eq:c:XaXb}}{=} \vkap(\alpha^\sharp,\alpha^\sharp) = \alpha(\alpha^\sharp)$\,; because of \,$\alpha(\alpha^\sharp) = \|\alpha\|^2 \neq 0$\,, the
statement follows. \emph{For (c),} we note that because the complex conjugation
\,$X \mapsto \overline{X}$\, is an involutive Lie algebra automorphism of \,$\lieg^{\C}$\, which leaves
\,$\liet$\, invariant and satisfies \,$\overline{\alpha^\sharp} = -\alpha^\sharp$\,, we have \,$\overline{X_\alpha} \in \lieg_{-\alpha}^{\C}$\, and therefore
there exists \,$a \in \C^\times$\, so that \,$\overline{X_\alpha} = a\cdot X_{-\alpha}$\, holds. We have \,$1 \overset{(b)}{=} \vkap(X_\alpha,X_{-\alpha}) = 
\tfrac{1}{a} \vkap(X_\alpha,\overline{X_\alpha})$\, and therefore \,$a = \vkap(X_\alpha,\overline{X_\alpha})$\,. We have \,$0 \neq X_\alpha + \overline{X_{\alpha}} \in \lieg$\,
and therefore, because \,$\lieg$\, is compact and hence \,$\vkap$\, is negative definite on \,$\lieg$\,, \,$0> \vkap(X_\alpha+\overline{X_{\alpha}},X_\alpha+\overline{X_{\alpha}})
= 2\,\vkap(X_\alpha,\overline{X_\alpha}) = 2a$\,, thus \,$a < 0$\,. \emph{For (d),} we have by (c) \,$V_\alpha(c) \in (\lieg_\alpha^{\C} \oplus \lieg_{-\alpha}^{\C})
\cap \lieg = \lieg_\alpha$\,, and therefore \,$\Menge{V_\alpha(c)}{c \in \C} \subset \lieg_\alpha$\,; because \,$\lieg_\alpha$\, is real-2-dimensional, equality follows.
\emph{For (e) and (f).} These follow essentially from the Jacobi identity, see \cite{Helgason:1978}: Lemma~III.5.1, p.~171 and Lemma~III.5.3, p.~172.
\emph{For (g).} See \cite{Knapp:2002}, Theorem~6.6, p.~351.
\beweisende

We next describe the transformation behavior of the Chevalley bases and the corresponding Chevalley constants:

\begin{Prop}
\label{P:c:chevtrafo}
Let \,$(X_\alpha)$\, be a Chevalley basis of \,$\lieg^{\C}$\, with the corresponding Chevalley constants \,$(c_{\alpha,\beta})$\,.
\begin{enumerate}
\item
Let constants \,$z_\alpha \in \C^\times$\, for every \,$\alpha \in \Delta$\, be given, so that the following properties are satisfied:
\begin{itemize}
\item[(i)] For every \,$\alpha \in \Delta$\, we have \,$z_\alpha \cdot z_{-\alpha} = 1$\,.
\item[(ii)] For every \,$\alpha,\beta \in \Delta$\, with \,$\alpha+\beta\in\Delta$\, we have \,$\eps_{\alpha,\beta} := \tfrac{z_\alpha \cdot z_\beta}{z_{\alpha+\beta}} \in \{\pm 1\}$\,.
\end{itemize}
Then \,$(z_\alpha\cdot X_\alpha)_{\alpha\in \Delta}$\, is another Chevalley basis of \,$\lieg^{\C}$\,, the corresponding Chevalley constants are
\,$(\eps_{\alpha,\beta}\cdot c_{\alpha,\beta})_{\alpha,\beta \in \Delta}$\,. 
\item
Every Chevalley basis of \,$\lieg^{\C}$\, is obtained by the construction of (a).
\end{enumerate}
\end{Prop}

\beweis
\emph{For (a).}
Put \,$\wt{X}_\alpha := z_\alpha \cdot X_\alpha$\, and \,$\wt{c}_{\alpha,\beta} := \eps_{\alpha,\beta}\cdot c_{\alpha,\beta}$\, for \,$\alpha,\beta\in\Delta$\,.
Obviously \,$\wt{X}_\alpha \in \lieg_\alpha^{\C}$\, holds for all \,$\alpha\in\Delta$\,, the numbers \,$\wt{c}_{\alpha,\beta}$\, are real by property (ii), and
it easily follows from (i) and (ii) that \,$((\wt{X}_\alpha),(\wt{c}_{\alpha,\beta}))$\, satisfies Equation~\eqref{eq:c:XaXb}. Moreover, for any \,$\alpha,\beta \in \Delta$\,
with \,$\alpha+\beta\in \Delta$\, we have
$$ \eps_{-\alpha,-\beta} = \tfrac{z_{-\alpha}\cdot z_{-\beta}}{z_{-(\alpha+\beta)}} \overset{\textrm{(i)}}{=} \tfrac{z_{\alpha+\beta}}{z_\alpha\cdot z_\beta}
= \eps_{\alpha,\beta}^{-1} \overset{\textrm{(ii)}}{=} \eps_{\alpha,\beta} \; ; $$
Therefrom it follows that \,$((\wt{X}_\alpha),(\wt{c}_{\alpha,\beta}))$\, also satisfies Equation~\eqref{eq:c:c-a-b}.

\emph{For (b).} 
Let two Chevalley bases \,$(X_\alpha)$\, and \,$(\wt{X}_\alpha)$\, of \,$\lieg^{\C}$\, (with corresponding Chevalley constants \,$(c_{\alpha,\beta})$\,
resp.~\,$(\wt{c}_{\alpha,\beta})$\,) be given. For each \,$\alpha \in \Delta$\,, the non-zero vectors \,$X_\alpha$\, and \,$\wt{X}_\alpha$\, lie in the complex-1-dimensional
root space \,$\lieg_\alpha^{\C}$\,, so there exists \,$z_\alpha \in \C^\times$\, so that \,$\wt{X}_\alpha = z_\alpha\cdot X_\alpha$\, holds.
It remains to show that the constants \,$(z_\alpha)$\, satisfy the conditions (i) and (ii) of (a). For (i): For any \,$\alpha \in \Delta$\, we have
$$ \alpha^\sharp \overset{\eqref{eq:c:XaXb}}{=} [\wt{X}_\alpha,\wt{X}_{-\alpha}] = z_\alpha\,z_{-\alpha}\,[X_\alpha,X_{-\alpha}] 
\overset{\eqref{eq:c:XaXb}}{=} z_\alpha\,z_{-\alpha}\,\alpha^\sharp $$
and therefore \,$z_\alpha\cdot z_{-\alpha} = 1$\,. For (ii): For any \,$\alpha,\beta \in \Delta$\, with \,$\alpha+\beta\in \Delta$\, we have
$$ X_{\alpha+\beta} = \tfrac{1}{z_{\alpha+\beta}} \,\wt{X}_{\alpha+\beta} 
\overset{\eqref{eq:c:XaXb}}{=} \tfrac{1}{z_{\alpha+\beta}}\,\tfrac{1}{\wt{c}_{\alpha,\beta}}\,[\wt{X}_\alpha,\wt{X}_\beta] 
= \tfrac{z_\alpha\cdot z_\beta}{z_{\alpha+\beta}} \, \tfrac{1}{\wt{c}_{\alpha,\beta}}\,[X_\alpha,X_\beta] 
\overset{\eqref{eq:c:XaXb}}{=} \tfrac{z_\alpha\cdot z_\beta}{z_{\alpha+\beta}} \, \tfrac{c_{\alpha,\beta}}{\wt{c}_{\alpha,\beta}}\,X_{\alpha+\beta} $$
and therefore \,$\tfrac{z_\alpha\cdot z_\beta}{z_{\alpha+\beta}} = \tfrac{\wt{c}_{\alpha,\beta}}{c_{\alpha,\beta}}$\,. It is a consequence of 
Proposition~\ref{P:c:c}(g) that \,$\wt{c}_{\alpha,\beta} = \pm c_{\alpha,\beta}$\, holds, and therefore we have 
\,$\tfrac{z_\alpha\cdot z_\beta}{z_{\alpha+\beta}} \in \{\pm 1\}$\,. 
\beweisende

The following proposition describes a way to choose a Chevalley basis in such a way that it is adapted to the position of the compact Lie algebra
\,$\lieg$\, within \,$\lieg^{\C}$\,, see property (i) in the proposition.

\begin{Prop}
\label{P:c:choice}
For every non-simple, positive root \,$\alpha \in \Delta_+ \setminus \Pi$\,, fix a decomposition \,$\alpha = \zeta_\alpha + \eta_\alpha$\, with \,$\zeta_\alpha,
\eta_\alpha \in \Delta_+$\,. Then there exists a Chevalley basis \,$(X_\alpha)$\, (with corresponding Chevalley constants 
\,$(c_{\alpha,\beta})$\,) with the following properties:
\begin{itemize}
\item[(i)]
For every \,$\alpha \in \Delta_+$\, we have \,$X_{-\alpha} = -\overline{X_\alpha}$\,. (Compare Proposition~\ref{P:c:c}(c).)
\item[(ii)]
For every \,$\alpha \in \Delta_+ \setminus \Pi$\, we have \,$c_{\zeta_\alpha,\eta_\alpha} > 0$\,. 
\end{itemize}
Any two such Chevalley bases have the same Chevalley constants \,$(c_{\alpha,\beta})$\,. 
\end{Prop}

\begin{Remark}
\label{R:c:decomp}
For each root \,$\alpha \in \Delta_+ \setminus \Pi$\, a decomposition \,$\alpha=\zeta_\alpha + \eta_\alpha$\, as required in Proposition~\ref{P:c:choice}
indeed exists: For example, it is well-known that for each \,$\alpha$\, there exists a simple root \,$\eta_\alpha$\, such that \,$\zeta_\alpha
:= \alpha - \eta_\alpha$\, is again a positive root.
\end{Remark}

{
\footnotesize
\emph{Proof of Proposition~\ref{P:c:choice}.}
We will show by induction on \,$\ell \geq 1$\, that there exists a Chevalley basis \,$(X_\alpha)$\, of \,$\lieg^{\C}$\,, 
with corresponding Chevalley constants \,$(c_{\alpha,\beta})$\,,
which has the property (i), and also satisfies \,$c_{\zeta_\alpha,\eta_\alpha} > 0$\, for all \,$\alpha \in \Delta_+$\, with \,$2 \leq \ell(\alpha) \leq \ell$\,.

First suppose \,$\ell=1$\,. Then we are to show that there exists a Chevalley basis of \,$\lieg^{\C}$\, which satisfies (i). For this we let an
arbitrary Chevalley basis \,$(X_\alpha)$\, of \,$\lieg^{\C}$\, be given, and denote the corresponding Chevalley constants by \,$(c_{\alpha,\beta})$\,. 
By Proposition~\ref{P:c:c}(c) there exists for every \,$\alpha \in \Delta$\, some \,$a_\alpha < 0$\, with \,$\overline{X_\alpha} = a_\alpha \cdot X_{-\alpha}$\,. 
We have \,$X_\alpha = \overline{\overline{X_\alpha}} = a_\alpha\,a_{-\alpha} X_\alpha$\, and therefore 
\begin{equation}
\label{eq:c:choice:a-alpha--alpha}
a_{\alpha} \cdot a_{-\alpha} = 1 \; . 
\end{equation}
Next we will show that for any \,$\alpha,\beta\in\Delta$\, with \,$\alpha+\beta\in\Delta$\, we have
\begin{equation}
\label{eq:c:choice:a-alpha-beta}
a_{\alpha+\beta} = -a_\alpha\cdot a_\beta \, . 
\end{equation}
For this we calculate \,$[\overline{X_\alpha},[X_\alpha,X_\beta]]$\, in two different ways: On one hand we have
$$ [\overline{X_\alpha},[X_\alpha,X_\beta]] = [a_\alpha \, X_{-\alpha}, [X_\alpha,X_\beta]] \overset{\eqref{eq:c:XaXb}}{=} a_\alpha \, c_{\alpha,\beta} [X_{-\alpha},X_{\alpha+\beta}]
\overset{\eqref{eq:c:XaXb}}{=} a_\alpha\, c_{\alpha,\beta} \, c_{-\alpha,(\alpha+\beta)} \, X_\beta \overset{(*)}{=} a_\alpha \, c_{\alpha,\beta}^2 \, X_\beta \; , $$
for the equals sign marked $(*)$ notice that we have \,$-\alpha + (\alpha+\beta) -\beta = 0$\, and therefore by Proposition~\ref{P:c:c}(e),(a) and
Equation~\eqref{eq:c:c-a-b}: \,$c_{-\alpha,(\alpha+\beta)} = c_{-\beta,-\alpha} = -c_{-\alpha,-\beta} = c_{\alpha,\beta}$\,. On the other hand we also have
\begin{align*}
[\overline{X_\alpha},[X_\alpha,X_\beta]] 
& = \overline{[X_\alpha, \overline{[X_\alpha,X_\beta]}]} 
\overset{\eqref{eq:c:XaXb}}{=} c_{\alpha,\beta} \, \overline{[X_\alpha, \overline{X_{\alpha+\beta}}]} 
= c_{\alpha,\beta} \, a_{\alpha+\beta} \, \overline{[X_\alpha, X_{-\alpha-\beta}]} 
\overset{\eqref{eq:c:XaXb}}{=} c_{\alpha,\beta} \, a_{\alpha+\beta} \, c_{\alpha,(-\alpha-\beta)} \, \overline{X_{-\beta}} \\
& \overset{(\dagger)}{=} -c_{\alpha,\beta}^2 \, a_{\alpha+\beta} \, a_{-\beta} \, X_{\beta}
\overset{\eqref{eq:c:choice:a-alpha--alpha}}{=} - \tfrac{a_{\alpha+\beta}}{a_\beta} c_{\alpha,\beta}^2 \, X_\beta \; ,
\end{align*}
for the equal signs marked \,$(\dagger)$\, notice that because of \,$\alpha + (-\alpha-\beta) + \beta=0$\, we have by Proposition~\ref{P:c:c}(e),(a)
\,$c_{\alpha,(-\alpha-\beta)} = c_{\beta,\alpha} = -c_{\alpha,\beta}$\,. By comparing the preceding two calculations we obtain Equation~\eqref{eq:c:choice:a-alpha-beta}.

Now put \,$z_\alpha := 1/\sqrt{-a_\alpha}$\, for every \,$\alpha \in \Delta$\,. Then it follows from Equation~\eqref{eq:c:choice:a-alpha--alpha} that we 
have \,$z_\alpha\cdot z_{-\alpha} = 1$\,, and it follows from Equation~\eqref{eq:c:choice:a-alpha-beta} that we have 
\,$\left(\tfrac{z_\alpha\cdot z_\beta}{z_{\alpha+\beta}}\right)^2=1$\,, hence \,$\tfrac{z_\alpha\cdot z_\beta}{z_{\alpha+\beta}} \in \{\pm 1\}$\,. 
By Proposition~\ref{P:c:chevtrafo}(a), \,$(\wt{X}_\alpha)$\, with \,$\wt{X}_\alpha := z_\alpha\cdot X_\alpha$\, therefore is another Chevalley basis of \,$\lieg^{\C}$\,,
and we have 
$$ \overline{\wt{X}_\alpha} = \overline{z_\alpha\,X_\alpha} = z_\alpha\,a_\alpha\,X_{-\alpha} = -\tfrac{1}{z_\alpha}\,X_{-\alpha} 
= -z_{-\alpha}\,X_{-\alpha} = -\wt{X}_{-\alpha} \;, $$
hence the new Chevalley basis \,$(\wt{X}_\alpha)$\, enjoys property (i) of the proposition.

We now suppose that \,$(X_\alpha)$\, is a Chevalley basis of \,$\lieg^{\C}$\,, with Chevalley constants \,$(c_{\alpha,\beta})$\,,
which satisfies property (i), and further satisfies \,$c_{\zeta_\alpha,\eta_\alpha} > 0$\, for every \,$\alpha\in \Delta_+$\, with
\,$2 \leq \ell(\alpha) \leq \ell-1$\,. Then put
$$ z_\alpha := \begin{cases}
1 & \text{for \,$\alpha\in\Delta_+$\, with \,$\ell(\alpha)\neq\ell$\,} \\
\sign(c_{\zeta_\alpha,\eta_\alpha}) & \text{for \,$\alpha\in\Delta_+$\, with \,$\ell(\alpha)=\ell$\,} \\
z_{-\alpha} & \text{for \,$\alpha \in -\Delta_+$\,}
\end{cases} \; . $$
In this way we have \,$z_\alpha \in \{\pm 1\}$\, in any case, and therefore \,$z_\alpha \cdot z_{-\alpha}=1$\, and \,$\eps_{\alpha,\beta} := 
\tfrac{z_\alpha\cdot z_\beta}{z_{\alpha+\beta}} \in \{\pm 1\}$\,. Therefore Proposition~\ref{P:c:chevtrafo}(a) shows that with \,$\wt{X}_\alpha := z_\alpha\cdot X_\alpha$\,
and \,$\wt{c}_{\alpha,\beta} := \eps_{\alpha,\beta}\cdot c_{\alpha,\beta}$\,, \,$(\wt{X}_\alpha)$\, is another Chevalley basis of \,$\lieg^{\C}$\,, 
its Chevalley constants being \,$(\wt{c}_{\alpha,\beta})$\,. For any \,$\alpha \in \Delta$\, we have \,$\wt{X}_{-\alpha} = z_{-\alpha} \cdot X_{-\alpha}
= -z_\alpha\cdot \overline{X_\alpha} = -\overline{\wt{X}_\alpha}$\,, hence the Chevalley basis \,$(\wt{X}_\alpha)$\, satisfies property (i). Moreover,
for any \,$\alpha \in \Delta$\, with \,$2 \leq \ell(\alpha) \leq \ell-1$\, we have \,$z_\alpha=z_{\zeta_\alpha} = z_{\eta_\alpha} = 1$\,, hence 
 \,$\wt{c}_{\zeta_\alpha,\eta_\alpha} = \eps_{\zeta_\alpha,\eta_\alpha} \, c_{\zeta_\alpha,\eta_\alpha} = c_{\zeta_\alpha,\eta_\alpha} > 0$\,, and for any \,$\alpha \in \Delta$\,
with \,$\ell(\alpha)=\ell$\, we have \,$z_{\zeta_\alpha}=z_{\eta_\alpha}=1$\, and \,$z_\alpha = \sign(c_{\zeta_\alpha,\eta_\alpha})$\,, hence
\,$\eps_{\zeta_\alpha,\eta_\alpha} = \sign(c_{\zeta_\alpha,\eta_\alpha})$\, and therefore \,$\wt{c}_{\zeta_\alpha,\eta_\alpha} = \eps_{\zeta_\alpha,\eta_\alpha} \,
c_{\zeta_\alpha,\eta_\alpha} > 0$\,. Thus we have shown \,$c_{\zeta_\alpha,\eta_\alpha} > 0$\, for all \,$\alpha \in \Delta$\, with 
\,$2 \leq \ell(\alpha) \leq \ell$\,. 

We postpone the proof of the uniqueness statement for the Chevalley constants until after we have given the algorithm for calculating these constants.
\hfill $\Box$

}

\medskip

To completely describe the Lie bracket of \,$\lieg^{\C}$\,, it suffices to fix a Chevalley basis \,$(X_\alpha)$\, of \,$\lieg^{\C}$\,, and calculate the corresponding 
Chevalley constants \,$(c_{\alpha,\beta})$\,. Up to sign, they are determined by Proposition~\ref{P:c:c}(g).  
However, to determine the sign correctly, we have to invest a bit more work. We now show how to do it, if \,$(X_\alpha)$\, is a Chevalley basis
of the kind of Proposition~\ref{P:c:choice}, corresponding to a family of decompositions \,$(\alpha = \zeta_\alpha + \eta_\alpha)_{\alpha \in \Delta_+\setminus\Pi}$\,.

By Proposition~\ref{P:c:c}(e),(a) and Equation~\eqref{eq:c:c-a-b} we have for any \,$\alpha,\beta\in\Delta_+$\, with \,$\beta\neq\alpha$\,
\begin{equation}
\label{eq:c:neg-reduction}
c_{\alpha,-\beta} = -c_{-\alpha,\beta} = \begin{cases}
c_{(\beta-\alpha),\alpha} & \text{if \,$\beta-\alpha \in \Delta_+$\,} \\
c_{(\alpha-\beta),\beta} & \text{if \,$\beta-\alpha \in (-\Delta_+)$\,} \\
0 & \text{otherwise}
\end{cases} 
\qmq{and} c_{-\alpha,-\beta} = -c_{\alpha,\beta} \; . 
\end{equation}
Therefore it suffices to calculate \,$c_{\alpha,\beta}$\, for \,$\alpha,\beta\in\Delta_+$\,. This is achieved by the following algorithm.
\smallskip
\begin{enumerate}
\item[\textbf{(C1)}] [\,$c_{\alpha,\beta}$\, with \,$\alpha+\beta \not\in \Delta$\,.] Iterate the following for all \,$\alpha,\beta \in \Delta_+$\,:
If \,$\alpha+\beta \not\in \Delta_+$\,, put \,$c_{\alpha,\beta} := 0$\,.
\item[\textbf{(C2)}] [Iterate on level.] Iterate steps (C3)--(C8) for \,$\ell=2,\dotsc,L$\,, where \,$L$\, denotes the maximal level of roots occurring in \,$\Delta$\,.
\item[\textbf{(C3)}] [Iterate on roots of level \,$\ell$\,.] Iterate steps (C4)--(C8) with \,$\alpha$\, running through all positive roots in \,$\Delta$\,
of level \,$\ell$\,.
\item[\textbf{(C4)}] 
Set \,$\zeta := \zeta_\alpha$\, and \,$\eta := \eta_\alpha$\,. 
\item[\textbf{(C5)}] [Calculate \,$c_{\zeta,\eta}$\, and \,$c_{\eta,\zeta}$\,.] Let \,$p$\, be the smallest integer so that \,$\eta-(p+1)\,\zeta \not\in \Delta$\, holds,
and put \,$q := p-2 \tfrac{\g{\eta}{\zeta}}{\|\zeta\|^2}$\,,
$$ c_{\zeta,\eta} := \sqrt{ \frac{q\cdot(1+p)}{2}} \cdot \|\zeta\|  $$
and \,$c_{\eta,\zeta} := -c_{\zeta,\eta}$\,. 
\item[\textbf{(C6)}] [Iterate on the decompositions of \,$\lambda$\,.] Iterate step (C7) for all pairs \,$(\gamma,\delta)$\, of positive roots
with \,$\gamma+\delta = \alpha$\, and \,$\gamma,\delta \not\in \{\zeta,\eta\}$\,.
\item[\textbf{(C7)}] [Calculate \,$c_{\gamma,\delta}$\,.]
Put 
{\scriptsize
\begin{gather*}
c_{11} := \begin{cases}
c_{\eta,(\gamma-\eta)} & \text{if \,$\gamma-\eta \in \Delta_+$\,} \\
c_{\gamma,(\eta-\gamma)} & \text{if \,$\gamma-\eta \in -\Delta_+$\,} \\
0 & \text{otherwise}
\end{cases} \; , \quad
c_{12} := \begin{cases}
c_{\delta,(\zeta-\delta)} & \text{if \,$\zeta-\delta \in \Delta_+$\,} \\
c_{\zeta,(\delta-\zeta)} & \text{if \,$\zeta-\delta \in -\Delta_+$\,} \\
0 & \text{otherwise}
\end{cases} \; , \\
c_{21} := \begin{cases}
c_{(\zeta-\gamma),\gamma} & \text{if \,$\zeta-\gamma \in \Delta_+$\,} \\
c_{(\gamma-\zeta),\zeta} & \text{if \,$\zeta-\gamma \in -\Delta_+$\,} \\
0 & \text{otherwise}
\end{cases}
\qmq{and}
c_{22} := \begin{cases}
c_{\eta,(\delta-\eta)} & \text{if \,$\delta-\eta \in \Delta_+$\,} \\
c_{\delta,(\eta-\delta)} & \text{if \,$\delta-\eta \in -\Delta_+$\,} \\
0 & \text{otherwise}
\end{cases}
\; . 
\end{gather*}

}

Then put \,$c_{\gamma,\delta} := \tfrac{1}{c_{\zeta,\eta}} \cdot (c_{11}\,c_{12} + c_{21}\,c_{22})$\,. 
\item[\textbf{(C8)}] (End of loops.)
\end{enumerate}

\bigskip

{\footnotesize
\emph{Proof for the correctness of the algorithm.} 
It is clear that the assignment \,$c_{\alpha,\beta}:=0$\, for every \,$\alpha,\beta\in\Delta_+$\, with \,$\alpha+\beta\not\in\Delta_+$\, in step (C1) is correct,
and that within the loop (C2)--(C8), every \,$c_{\alpha,\beta}$\, with \,$\alpha,\beta,\alpha+\beta\in\Delta_+$\, is assigned to exactly once (either in step (C5)
or in step (C7)). We have to show that the latter assignments are correct.

In step (C5), the numbers \,$p$\, and \,$q$\, are chosen such that \,$\Menge{\eta+k\,\zeta}{-p\leq k\leq q}$\, is the \,$\zeta$-string through \,$\eta$\,
(for the correctness of \,$q$\,, see Equation~\eqref{eq:posrts:stringpq}). Therefore, and because we have \,$c_{\zeta,\eta}>0$\, by definition,
Proposition~\ref{P:c:c}(g) shows that the assignment to \,$c_{\zeta,\eta}$\, in step (C5) is correct.
The correctness of the following assignment \,$c_{\eta,\zeta}:=-c_{\zeta,\eta}$\, follows by Proposition~\ref{P:c:c}(a).

For the correctness of the assignment in step (C7): In the situation of that step, we have \,$\zeta+\eta = \alpha = \gamma+\delta$\, and therefore
\,$-\zeta-\eta+\gamma+\delta=0$\,, and none of the four summands in the latter equation is the negative of one of the others. Therefore Proposition~\ref{P:c:c}(f)
shows that we have
$$ c_{-\zeta,-\eta}\,c_{\gamma,\delta} + c_{-\eta,\gamma}\,c_{-\zeta,\delta} + c_{\gamma,-\zeta}\,c_{-\eta,\delta} = 0 $$
and thus  (note that \,$c_{-\zeta,-\eta}=-c_{\zeta,\eta}$\, by Equation~\eqref{eq:c:c-a-b}) 
\begin{equation}
\label{eq:c:algoproof-cgammadelta}
c_{\gamma,\delta} = \frac{1}{c_{\zeta,\eta}} \cdot (c_{-\eta,\gamma}\,c_{-\zeta,\delta} + c_{\gamma,-\zeta}\,c_{-\eta,\delta}) \; . 
\end{equation}

If \,$\gamma-\eta$\, is not a root, then we have \,$c_{-\eta,\gamma}=0$\,. Otherwise we have 
either \,$\gamma-\eta \in \Delta_+$\,, or else \,$\gamma-\eta \in -\Delta_+$\, and then \,$\eta-\gamma \in \Delta_+$\,.
In either case, by Proposition~\ref{P:c:c}(e),(a) and Equation~\eqref{eq:c:c-a-b}
\begin{equation}
\label{eq:c:algoproof-cetagamma}
c_{-\eta,\gamma} \overset{(e)}{=} c_{\gamma,(\eta-\gamma)} \overset{\eqref{eq:c:c-a-b}}{=} -c_{-\gamma,(\gamma-\eta)} \overset{(a)}{=} c_{(\gamma-\eta),-\gamma}
\overset{(e)}{=} c_{\eta,(\gamma-\eta)}
\end{equation}
holds. Notice that the level of \,$\eta+(\gamma-\eta)=\gamma$\, (if \,$\gamma-\eta$\, is positive) or of \,$\gamma+(\eta-\gamma)=\eta$\, 
(if \,$\gamma-\eta$\, is negative) is strictly less than \,$\ell$\,, and therefore the value of \,$c_{\eta,(\gamma-\eta)}$\, resp.~of \,$c_{\gamma,(\eta-\gamma)}$\, 
has already been calculated by the algorithm in a previous iteration. Therefore Equation~\eqref{eq:c:algoproof-cetagamma} shows that the value \,$c_{11}$\,
computed in step (C7) equals \,$c_{-\eta,\gamma}$\,. 

Analogously one sees that \,$c_{12}=c_{-\zeta,\delta}$\,,\,$c_{21} = c_{\gamma,-\zeta}$\, and \,$c_{22} = c_{-\eta,\delta}$\, holds. Therefore
Equation~\eqref{eq:c:algoproof-cgammadelta} shows that the value of \,$c_{\gamma,\delta}$\, is calculated correctly in step (C7) of the algorithm.
\strut\hfill $\Box$

}

\bigskip

{
\footnotesize
\emph{Proof of the uniqueness statement for the Chevalley constants in Proposition~\ref{P:c:choice}.}
This is also a consequence of the proof of the correctness of the algorithm. Indeed that proof shows that the Chevalley constants \,$(c_{\gamma,\delta})$\, are determined
uniquely by the decompositions \,$(\alpha = \zeta_\alpha+\eta_\alpha)$\,, as was claimed in Proposition~\ref{P:c:choice}.
\strut\hfill $\Box$

}

\section{Reconstructing the symmetric involution}
\label{Se:sigma}

We now suppose that a Riemannian symmetric space \,$M=G/K$\, of compact type is given. Then \,$G$\, is a semisimple Lie group, and the symmetric
structure of \,$M$\, is given by an involutive Lie algebra automorphism \,$\sigma$\, of the Lie algebra \,$\lieg$\, of \,$G$\,. It will now be our objective
to describe how to reconstruct the action of \,$\sigma$\, on \,$\lieg$\, from the information contained in the Satake diagram of \,$M$\,. 

\,$\sigma$\, induces the splitting \,$\lieg = \liek\oplus \liem$\, of \,$\lieg$\,, where \,$\liek := \Eig(\sigma,1)$\, is the Lie algebra of the 
isotropy group \,$K$\, and \,$\liem := \Eig(\sigma,-1)$\, is a linear subspace of \,$\lieg$\, which is isomorphic to the tangent space \,$T_pM$\,
in a canonical way.

We let \,$\liea$\, be a maximal flat subspace of \,$\liem$\, and let 
\,$\liet$\, be a Cartan subalgebra of \,$\lieg$\, with \,$\liea \subset \liet$\,. Then \,$\liet$\, is invariant under
\,$\sigma$\, (see \cite{Loos:1969-2}, Lemma~VI.3.2, p.~72). 

In the sequel, we apply the results of the preceding sections in this situation. For this we again also consider the complexification \,$\lieg^{\C}$\,
of \,$\lieg$\, as a complex Lie algebra. We denote the complexification of \,$\sigma$\, also by \,$\sigma$\,; this is
an involutive Lie algebra automorphism of \,$\lieg^{\C}$\,. 

The action of \,$\sigma$\, on \,$\liet$\, has been described by \textsc{Satake} in the following way
(see, for example, \cite{Loos:1969-2}, Section~VII.3.3, p.~132ff.):
There exists a partition of the set \,$\Pi$\, of simple roots of \,$\lieg$\, into two subsets:
\,$\Pi= \{\alpha_1,\dotsc,\alpha_r\} \dot\cup \{\beta_1,\dotsc,\beta_s\}$\, (with \,$r,s \in \N_0$\,, \,$r+s=\rk(\lieg)$\,) and an involutive permutation
\,$\pi \in \mathfrak{S}_r$\, so that for each \,$\alpha_k$\, (\,$k=1,\dotsc,r$\,) we have
\begin{equation}
\label{eq:sigma:sigmaalpha}
\sigma(\alpha_k^\sharp) = -\alpha_{\pi(k)}^\sharp - \sum_{j =1}^s n_{kj}\,\beta_j^\sharp
\end{equation}
with some non-negative integers \,$n_{kj}$\,, whereas for each \,$\beta_k$\, (\,$k=1,\dotsc,s$\,) we have
\begin{equation}
\label{eq:sigma:sigmabeta}
\sigma(\beta_k^\sharp) = \beta_k^\sharp \; .
\end{equation}

The \emph{Satake diagram} of the symmetric space \,$M$\, is obtained by ``annotating'' the Dynkin diagram of \,$\lieg$\, in the following way: We color the node
corresponding to each simple root of \,$\lieg$\, either white or black, according to whether it is one of the \,$\alpha_k$\, or one of the \,$\beta_k$\,
(these simple roots are thus called \emph{white roots} and \emph{black roots}, respectively), 
and wherever two white roots are interchanged by the involutive permutation \,$\pi$\,,
we indicate this fact by drawing a curved both-sided arrow between them. 

\medskip

We now show how to reconstruct the action of \,$\sigma$\, on \,$\lieg^{\C}$\, (or on \,$\lieg$\,) from the Satake diagram of \,$M$\,. 

First, it was already described by Satake how to reconstruct the action of \,$\sigma$\, on \,$\liet$\,: For this, we note that the partition 
\,$\Pi= \{\alpha_1,\dotsc,\alpha_r\} \dot\cup \{\beta_1,\dotsc,\beta_s\}$\, of the set of simple roots into the sets of white resp.~black roots,
and the involutive permutation \,$\pi \in \mathfrak{S}_r$\, can be read off the Satake diagram of \,$M$\, immediately. To apply formulas~\eqref{eq:sigma:sigmaalpha}
and \eqref{eq:sigma:sigmabeta} to obtain the action of \,$\sigma$\, on the simple roots, we therefore only need to determine the numbers \,$n_{kj}$\, 
occurring in Equation~\eqref{eq:sigma:sigmaalpha}, and this can be done in the following way:
From Equation~\eqref{eq:sigma:sigmaalpha} it follows, putting \,$\beta_\ell^* := \tfrac{2\,\beta_\ell^\sharp}{\|\beta_\ell^\sharp\|^2}$\,:
$$ \sum_{j=1}^s \beta_j(\beta_\ell^*)\,n_{kj} = -\alpha_k(\sigma(\beta_\ell^*)) - \alpha_{\pi(k)}(\beta_\ell^*) 
\overset{\eqref{eq:sigma:sigmabeta}}{=} -\alpha_k(\beta_\ell^*) - \alpha_{\pi(k)}(\beta_\ell^*) \; . $$
The numbers \,$\beta_j(\beta_\ell^*)$\,, \,$\alpha_k(\beta_\ell^*)$\, and \,$\alpha_{\pi(k)}(\beta_\ell^*)$\, are known from the Dynkin diagram of \,$\lieg$\,.
Because the matrix \,$(\beta_j(\beta_\ell^*))_{j,l=1,\dotsc,s}$\, is invertible (it is the Cartan matrix of
the Lie algebra \,$\liek^\liea = \Menge{X\in \liek}{[X,\liea]=0}$\,),
we can solve for \,$n_{kj}$\,. 
Because \,$\Pi$\, is a basis of \,$i\liet$\,, from Equations~\eqref{eq:sigma:sigmaalpha} and \eqref{eq:sigma:sigmabeta} 
we therefrom know the action of \,$\sigma$\, on all of \,$\liet^{\C}$\,.

But now we also want to describe the action of \,$\sigma$\, on the root spaces of \,$\lieg^{\C}$\,. For this, we let a Chevalley basis \,$(X_\alpha)$\, of \,$\lieg^{\C}$\,
of the kind of Proposition~\ref{P:c:choice} be given, in particular we have for any \,$\alpha\in\Delta$\,
\begin{equation}
\label{eq:sigma:Xconj}
X_{-\alpha} = -\overline{X_\alpha} \; .
\end{equation}
We denote the Chevalley constants corresponding to \,$(X_\alpha)$\, by \,$(c_{\alpha,\beta})$\,. Because \,$\sigma$\, is an involutive automorphism of \,$\lieg$\,,
we have for every \,$\alpha\in\Delta$\,: \,$\sigma(X_\alpha) \in \lieg^{\C}_{\alpha\circ \sigma^{-1}} = \lieg^{\C}_{\sigma(\alpha)}$\,
(where we write also \,$\sigma(\alpha)$\, for \,$\alpha\circ\sigma^{-1}$\, in the sequel; note that with this notation, \,$(\sigma(\alpha))^\sharp = \sigma(\alpha^\sharp)$\,
holds) and therefore there exists
\,$s_\alpha \in \C$\, with
\begin{equation}
\label{eq:sigma:sdef}
\sigma(X_\alpha) = s_\alpha\cdot X_{\sigma(\alpha)} \; .
\end{equation}
Clearly the action of \,$\sigma$\, on the root spaces is determined completely by the constants \,$s_\alpha$\,, and it is our objective in the following
to determine these constants.

\begin{Prop}
\label{P:sigma:s}
Suppose \,$\alpha,\beta\in\Delta$\,. 
\begin{enumerate}
\item \,$s_{-\alpha} = s_{\alpha}^{-1} = s_{\sigma(\alpha)}$\,.
\item \,$|s_\alpha|=1$\,.
\item If \,$\alpha+\beta\in\Delta$\, holds, we have \,$s_{\alpha+\beta} = \tfrac{c_{\sigma\alpha,\sigma\beta}}{c_{\alpha,\beta}}\,s_\alpha\,s_\beta$\,.
\item If \,$\sigma(\beta)=\beta$\, holds, we have \,$s_\beta=1$\,.
\end{enumerate}
\end{Prop}

\beweis
\emph{For (a).} By the involutivity of \,$\sigma$\, we have \,$X_\alpha = \sigma^2(X_\alpha) \overset{\eqref{eq:sigma:sdef}}{=}
\sigma(s_\alpha\,X_{\sigma(\alpha)}) \overset{\eqref{eq:sigma:sdef}}{=} s_\alpha\,s_{\sigma(\alpha)} \, X_\alpha$\, and therefore \,$1 = s_\alpha\,s_{\sigma(\alpha)}$\,. We also have
\,$\sigma(\alpha^{\sharp}) \overset{\eqref{eq:c:XaXb}}{=} \sigma([X_\alpha,X_{-\alpha}]) = [\sigma(X_\alpha),\sigma(X_{-\alpha})] \overset{\eqref{eq:sigma:sdef}}{=}
s_\alpha\,s_{-\alpha}\,[X_{\sigma(\alpha)},X_{-\sigma(\alpha)}] \overset{\eqref{eq:c:XaXb}}{=} s_\alpha\,s_{-\alpha}\,\sigma(\alpha^\sharp)$\, and hence
\,$1 = s_\alpha\,s_{-\alpha}$\,. 

\emph{For (b).} We have \,$\overline{s_\alpha}^{-1}\,X_{\sigma(\alpha)} \overset{(a)}{=} \overline{s_{-\alpha}}\,X_{\sigma(\alpha)}
\overset{\eqref{eq:sigma:Xconj}}{=} -\overline{s_{-\alpha}\,X_{-\sigma(\alpha)}} \overset{\eqref{eq:sigma:sdef}}{=} -\overline{\sigma(X_{-\alpha})}
= -\sigma(\overline{X_{-\alpha}}) \overset{\eqref{eq:sigma:Xconj}}{=} \sigma(X_\alpha) \overset{\eqref{eq:sigma:sdef}}{=} s_\alpha\,X_{\sigma(\alpha)}$\,
and therefore \,$\overline{s_\alpha}^{-1} = s_\alpha$\,, hence \,$|s_\alpha|=1$\,.  

\emph{For (c).} We have \,$s_{\alpha+\beta}\,X_{\sigma(\alpha+\beta)} \overset{\eqref{eq:sigma:sdef}}{=} \sigma(X_{\alpha+\beta}) \overset{\eqref{eq:c:XaXb}}{=}
\tfrac{1}{c_{\alpha,\beta}}\,\sigma([X_\alpha,X_\beta]) = \tfrac{1}{c_{\alpha,\beta}}\,[\sigma(X_\alpha),\sigma(X_\beta)] \overset{\eqref{eq:sigma:sdef}}{=}
\tfrac{1}{c_{\alpha,\beta}}\,s_\alpha\,s_\beta\,[X_{\sigma(\alpha)},X_{\sigma(\beta)}] \overset{\eqref{eq:c:XaXb}}{=} \tfrac{c_{\sigma\alpha,\sigma\beta}}{c_{\alpha,\beta}}\,
s_\alpha\,s_\beta\,X_{\sigma(\alpha)+\sigma(\beta)}$\, and therefore \,$s_{\alpha+\beta} = \tfrac{c_{\sigma\alpha,\sigma\beta}}{c_{\alpha,\beta}}\,s_\alpha\,s_\beta$\,.

\emph{For (d).} We now suppose that \,$\sigma(\beta)=\beta$\, holds. Then we have \,$s_\beta \in \{\pm 1\}$\, by (a).
Assume that \,$s_\beta = -1$\, holds. Then with \,$v := X_\beta - X_{-\beta} \overset{\eqref{eq:sigma:Xconj}}{\in} \lieg$\,
we have \,$\sigma(v) = -v$\, and therefore \,$v \in \liem$\,. By hypothesis we have \,$\beta\circ\sigma=\beta$\,, and therefore \,$\beta$\, vanishes on 
\,$\liea \subset \Eig(\sigma,-1)$\,. For every \,$H \in \liea$\, we therefore have \,$[H,v] = [H,X_\beta] - [H,X_{-\beta}] = \beta(H)\cdot(X_\beta + X_{-\beta}) = 0$\,,
and thus \,$\liea \oplus \R v$\, is a flat subspace of \,$\liem$\,, in contradiction to the maximality of \,$\liea$\,. Therefore we have \,$s_\beta=1$\,. 
\beweisende

\begin{Prop}
\label{P:sigma:alphaadap}
Suppose that \,$M$\, is irreducible.
\begin{enumerate}
\item 
\begin{enumerate}
\item If \,$M$\, is of one of the types \,$\textsf{AIII}(n,q) = \SU(n)/\mathrm{S}(\Ug(q)\times\Ug(n-q))$\, with \,$2q<n-1$\,, 
\,$\textsf{DIII}(n) = \SO(2n)/\Ug(n)$\, with \,$n$\, odd,
or \,$\textsf{EIII} =E_6/(U(1)\cdot \SO(10))$\,, we label simple roots as in the following Satake diagrams:

\begin{center}
\begin{longtable}{ccc}
\begin{minipage}{3cm}
\begin{center}
\,$\textsf{AIII}(n,q)$\, \\
with \,$2q<n-1$\, 
\end{center}
\end{minipage}
& 
\begin{minipage}{5cm}
\medskip
\xymatrix@=.4cm{
\displaystyle \mathop{\circ}^{\alpha_q}\ar@{-}[r]<-.8ex>\ar@{<->}@/^/[d]<-1ex>
        & \displaystyle \mathop{\circ}^{\alpha_{q-1}}\ar@{.}[r]<-.8ex>\ar@{<->}@/^/[d]<-1ex>
        & \displaystyle \mathop{\circ}^{\alpha_2}\ar@{-}[r]<-.8ex>\ar@{<->}@/^/[d]<-1ex>
        & \displaystyle \mathop{\circ}^{\alpha_1}\ar@{-}[r]<-.8ex>\ar@{<->}@/^/[d]<-1ex>
        & \displaystyle \mathop{\bullet}^{\beta_1} \ar@{.}[d] \\
{\displaystyle \mathop{\circ}_{\alpha_{\pi(q)}}}\ar@{-}[r]<1.2ex>
        & {\displaystyle \mathop{\circ}_{\alpha_{\pi(q-1)}}}\ar@{.}[r]<1.2ex>
        & {\displaystyle \mathop{\circ}_{\alpha_{\pi(2)}}}\ar@{-}[r]<1.2ex>
        & {\displaystyle \mathop{\circ}_{\alpha_{\pi(1)}}}\ar@{-}[r]<1.2ex>
        & {\displaystyle \mathop{\bullet}_{\beta_s}} \\
}
\medskip
\end{minipage}
& 
\begin{minipage}{4cm}
\begin{center}
\,$\beta := \beta_1+\dotsc+\beta_s \in \Delta_+$\, \\
\,$s = (n-1)-2q \geq 1$\,
\end{center}
\end{minipage}
\\

\begin{minipage}{3cm}
\begin{center}
\,$\textsf{DIII}(n)$\, \\
with \,$n$\, odd
\end{center}
\end{minipage}
&
\begin{minipage}{5cm}
\medskip
\xymatrix@=.4cm{
&&&&&&& \stackrel{\alpha_1}{\circ} \ar@{<->}@/^/[2,0] \\
\bullet \ar@{-}[r] & \circ \ar@{-}[r] & \bullet \ar@{-}[r] & \circ \ar@{.}[r] & \ldots \ar@{.}[r] & \circ \ar@{-}[r] & {\displaystyle \mathop{\bullet}^{\beta}_{\,}} \ar@{-}[ur] \ar@{-}[dr] & \\
&&&&&&& {\displaystyle \mathop{\circ}_{\alpha_{\pi(1)}}}
}
\medskip
\end{minipage}
& \\

\,$\textsf{EIII}$\, & 
\begin{minipage}{5cm}
\medskip
\xymatrix@=.4cm{
{\displaystyle \mathop{\circ}^{\alpha_1}} \ar@{-}[r]<-.8ex> \ar@{<->}@/^2pc/[0,4]
        & \displaystyle \mathop{\bullet}^{\beta_1} \ar@{-}[r]<-.8ex>
        & \displaystyle \mathop{\bullet}^{\beta_2} \ar@{-}[r]<-.8ex> \ar@{-}[d]
        & \displaystyle \mathop{\bullet}^{\beta_3} \ar@{-}[r]<-.8ex>
        & {\displaystyle \mathop{\circ}^{\alpha_{\pi(1)}}} \\
&& {\displaystyle \mathop{\circ}_{\alpha_{2}}} &&
}
\medskip
\end{minipage}
& \,$\beta := \beta_1 + \beta_2 + \beta_3 \in \Delta_+$\,
\end{longtable}
\end{center}
Then we have
$$ s_{\alpha_{\pi(1)}} = \frac{c_{\beta,\alpha_1}}{c_{\alpha_{\pi(1)},\beta}} \, s_{\alpha_{1}} \qmq{and} s_{\alpha_{\pi(k)}} = s_{\alpha_k} \text{\quad for \,$k\not\in \{1,\pi(1)\}$\,} \; . $$

\item
If \,$M$\, is of a type not mentioned in (i), we have \,$s_{\alpha_{\pi(k)}} = s_{\alpha_k}$\, for all \,$k$\,. 
\end{enumerate}
 
\item
\,$\sigma$\, is congruent to another involutive automorphism \,$\wt{\sigma}$\, of \,$\lieg$\, so that (with \,$\wt{s}_\alpha$\, being defined 
for \,$\wt{\sigma}$\, analogous to Equation~\eqref{eq:sigma:sdef}) we have:
\begin{enumerate}
\item If \,$M$\, is of one of the types \,$\textsf{AIII}(n,q)$\, with \,$2q<n-1$\,, \,$\textsf{DIII}(n)$\, with \,$n$\, odd, and \,\textsf{EIII}\,, then
\,$\wt{s}_{\alpha_{\pi(1)}} = \tfrac{c_{\beta,\alpha_1}}{c_{\alpha_{\pi(1)},\beta}}$\, and  \,$\wt{s}_{\alpha_k} = 1$\, for all \,$k \neq \pi(1)$\,
(where \,$\pi$\,, \,$\alpha_k$\, and \,$\beta$\, have the same meaning as in (a)(i)). 
\item If \,$M$\, is of any other type, then \,$\wt{s}_{\alpha_k}=1$\, for all \,$k$\,. 
\end{enumerate}
In either case, we have \,$\wt{s}_\alpha \in\{\pm 1\}$\, for all \,$\alpha\in\Delta$\,. 
\end{enumerate}
\end{Prop}

\begin{Remarks}
\begin{enumerate}
\item
If \,$M$\, is not irreducible, then the statement of Proposition~\ref{P:sigma:alphaadap} holds within each irreducible factor of \,$M$\,. 
\item
The quotient \,$\tfrac{c_{\beta,\alpha_1}}{c_{\alpha_{\pi(1)},\beta}}$\, occurring in Proposition~\ref{P:sigma:alphaadap} is always \,$\in \{\pm 1\}$\,, and the Chevalley
basis of \,$\lieg$\, can be chosen so that it equals \,$1$\,; such a Chevalley basis can, for example, be obtained by applying Proposition~\ref{P:c:choice} 
resp.~Algorithm~(C) with \,$\zeta_{\beta+\alpha_1} = \beta$\,, \,$\eta_{\beta+\alpha_1} = \alpha_1$\,, \,$\zeta_{\alpha_{\pi(1)}+\beta} = \alpha_{\pi(1)}$\, and
\,$\eta_{\alpha_{\pi(1)}+\beta} = \beta$\,. If we use such a Chevalley basis, then the involutive automorphism \,$\wt{\sigma}$\, of
Proposition~\ref{P:sigma:alphaadap}(b) satisfies \,$\wt{s}_{\alpha_k} = 1$\, for all \,$k$\, in any case.
\item
Even if the Chevalley basis is set up so that we have \,$\wt{s}_\alpha = 1$\, for every simple root \,$\alpha$\,, we do not generally have \,$\wt{s}_\alpha=1$\,
for all \,$\alpha\in\Delta$\,. 
\end{enumerate}
\end{Remarks}

\beweis
\emph{For (a).}
If \,$\pi = \id$\, holds (i.e.~the Satake diagram of \,$M$\, has no arrows), then there is nothing to show. If \,$M$\, has no black roots (this is also
the case if \,$M$\, is a simple Lie group seen as symmetric space), then it follows
from Equation~\eqref{eq:sigma:sigmaalpha} that for any \,$k \in \{1,\dotsc,r\}$\,, \,$\sigma(\alpha_k) = -\alpha_{\pi(k)}$\, holds, and therefore we have 
by Proposition~\ref{P:sigma:s}(a): \,$s_{\alpha_k} = s_{-\sigma(\alpha_k)} = s_{\alpha_{\pi(k)}}$\,. Thus it only remains to consider the spaces whose
Satake diagrams contain both arrows and black roots, and an inspection of the Satake diagrams of all the irreducible Riemannian symmetric spaces
(see for example \cite{Loos:1969-2}, p.~147f.) shows that they are only the spaces of type \,$\textsf{AIII}(n,q)$\, with \,$2q<n-1$\,, \,$\textsf{DIII}(n)$\, with \,$n$\, odd,
and \,$\textsf{EIII}$\,. 

Let us therefore now suppose that \,$M$\, is of one of these three types. Then we have \,$\sigma(\beta)=\beta$\,, \,$-\sigma(\alpha_1) = \alpha_{\pi(1)} + \beta$\, and
\,$-\sigma(\alpha_{\pi(1)}) = \alpha_1 + \beta$\, (with the roots \,$\alpha_1$\,, \,$\alpha_{\pi(1)}$\, and \,$\beta$\, being defined as in the relevant part of the proposition)
and therefore 
$$ s_{\alpha_1} 
\overset{(a)}{=} s_{-\sigma(\alpha_1)} 
= s_{\alpha_{\pi(1)} + \beta} 
\overset{(c)}{=} \frac{c_{\sigma(\alpha_{\pi(1)}),\sigma(\beta)}}{c_{\alpha_{\pi(1)},\beta}} \cdot s_{\alpha_{\pi(1)}} \cdot \underbrace{s_\beta}_{\underset{(d)}{=}1}
= \frac{c_{-(\alpha_1+\beta),\beta}}{c_{\alpha_{\pi(1)},\beta}} \cdot s_{\alpha_{\pi(1)}} \overset{(*)}{=} \frac{c_{\beta,\alpha_1}}{c_{\alpha_{\pi(1)},\beta}} \cdot s_{\alpha_{\pi(1)}} \; . 
$$
(Herein, the letters (a), (c) and (d) refer to the respective parts of Proposition~\ref{P:sigma:s}, and the equals sign marked $(*)$ follows from 
Proposition~\ref{P:c:c}(e).)

Moreover, if \,$M$\, is of the type \,$\textsf{AIII}(n,q)$\, with \,$2q<n-1$\,, then we have \,$\sigma(\alpha_k) = -\alpha_{\pi(k)}$\, for any \,$k \not\in \{1,\pi(1)\}$\,
and therefore for such \,$k$\, by Proposition~\ref{P:sigma:s}(a): \,$s_{\alpha_k} = s_{-\sigma(\alpha_k)} = s_{\alpha_{\pi(k)}}$\,. On the other hand, if \,$M$\,
is of type \,$\textsf{DIII}(2n+1)$\, or of type \,$\textsf{EIII}$\,, then we have \,$\pi(k)=k$\, for all \,$k \not\in \{1,\pi(1)\}$\,, and therefore
\,$s_{\alpha_{\pi(k)}} = s_{\alpha_k}$\, then trivially holds. 

\emph{For (b).} For arbitrary \,$H \in \liet$\,, the map \,$\ad(H)$\, is a derivation of \,$\lieg$\, with \,$\ad(H)|\liet=0$\,,
therefore \,$B := \exp(\ad(H))$\, is a Lie algebra automorphism of \,$\lieg$\, with \,$B|\liet = \id_\liet$\,.
Hence \,$\wt{\sigma} := B \circ \sigma \circ B^{-1}$\, then is another involutive automorphism of \,$\lieg$\, with \,$\wt{\sigma}|\liet = \sigma|\liet$\,. 
Thus \,$\wt{\sigma}$\, describes the same symmetric structure as \,$\sigma$\, does, and 
we can define \,$\wt{s}_\alpha$\, with respect
to \,$\wt{\sigma}$\, analogous to Equation~\eqref{eq:sigma:sdef}. Doing so, the results of the present section, especially part (a) of the present
proposition, are true mutatis mutandis with \,$\wt{s}_\alpha$\, in the place of \,$s_\alpha$\,. 

Because of (a) it therefore suffices to show that \,$H$\, 
can be chosen in such a way that for every \,$k \in \{1,\dotsc,r\}$\, at least one of the equations \,$\wt{s}_{\alpha_k}=1$\, and \,$\wt{s}_{\alpha_{\pi(k)}}=1$\, holds.
For this purpose, for every \,$k$\, we let \,$\tau(k)=\tau(\pi(k))$\, be an arbitrarily chosen element of \,$\{k,\pi(k)\}$\, and let \,$t_k \in \R$\, be such
that 
\begin{equation}
\label{eq:sigma:alphaadap:tk}
s_{\alpha_k} = e^{i\,t_k}
\end{equation}
holds. Then we let \,$H \in \liet$\, be the element characterized by
$$ \g{H}{\alpha_k^\sharp} = \g{H}{\alpha_{\pi(k)}^\sharp} = \tfrac12\,i\,t_{\tau(k)} \qmq{for every \,$k$\,} $$
and 
$$ \g{H}{\beta_k^\sharp} = 0 \qmq{for every black root \,$\beta_k$\,} \; . $$
Then we have for every \,$k \in \{1,\dotsc,r\}$\,
$$ \ad(H)\,X_{\alpha_k} = [H,X_{\alpha_k}] = \alpha_k(H)\,X_{\alpha_k} = \g{H}{\alpha_k^\sharp}\,X_{\alpha_k} = \tfrac12\,i\,t_{\tau(k)}\,X_{\alpha_k} $$
and for every \,$k\in\{1,\dotsc,s\}$\,
$$ \ad(H)\,X_{\beta_k} = [H,X_{\beta_k}] = \beta_k(H)\,X_{\beta_k} = \g{H}{\beta_k^\sharp} \, X_{\beta_k} = 0 \; . $$
Therefore we obtain
\begin{equation}
\label{eq:sigma:alphaadap:BXalpha}
B(X_{\alpha_k}) = e^{\tfrac12\,i\,t_{\tau(k)}} \cdot X_{\alpha_k} 
\end{equation}
and -- with the notation from Equation~\eqref{eq:sigma:sigmaalpha} --
\begin{equation}
\label{eq:sigma:alphaadap:BXsigmaalpha}
B(X_{\sigma(\alpha_k)}) = \g{H}{\sigma(\alpha_k^\sharp)}\,X_{\sigma(\alpha_k)} = \g{H}{-\alpha_{\pi(k)}^\sharp - \sum_{j=1}^s n_{kj}\,\beta_j^\sharp}\,X_{\sigma(\alpha_k)}
= e^{-\tfrac12\,i\,t_{\tau(k)}}\cdot X_{\sigma(\alpha_k)} \; . 
\end{equation}
Therefrom it follows that
\begin{align*}
\wt{\sigma}(X_{\alpha_{\tau(k)}}) 
= (B \circ \sigma \circ B^{-1})(X_{\alpha_{\tau(k)}}) 
& \overset{\eqref{eq:sigma:alphaadap:BXalpha}}{=} (B \circ \sigma)(e^{-\tfrac12\,i\,t_{\tau(k)}}\,X_{\alpha_{\tau(k)}})
\overset{\eqref{eq:sigma:sdef}}{=} B(s_{\alpha_{\tau(k)}} \, e^{-\tfrac12\,i\,t_{\tau(k)}}\,X_{\sigma(\alpha_{\tau(k)})}) \\
& \overset{\eqref{eq:sigma:alphaadap:BXsigmaalpha}}{=} e^{-i\,t_{\tau(k)}}\,s_{\alpha_{\tau(k)}}\,X_{\sigma(\alpha_{\tau(k)})}
\overset{\eqref{eq:sigma:alphaadap:tk}}{=} X_{\sigma(\alpha_{\tau(k)})} 
\end{align*}
and therefore \,$\wt{s}_{\alpha_{\tau(k)}} = 1$\, holds.
\beweisende

We are now able to reconstruct the action of \,$\sigma$\, on the root spaces of \,$\lieg^{\C}$\,. For this we may suppose without loss of generality
that the Chevalley basis \,$(X_\alpha)$\, and the involutive automorphism \,$\sigma$\, are adapted to each other in such a way that \,$\sigma$\,
equals the automorphism \,$\wt{\sigma}$\, from Proposition~\ref{P:sigma:alphaadap}(b). The following algorithm then computes the \,$s_\alpha$\, for
\,$\alpha \in \Delta_+$\,; by Proposition~\ref{P:sigma:s}(a) we then know \,$s_\alpha$\, for all \,$\alpha\in\Delta$\,.
\medskip
\begin{enumerate}
\item[\textbf{(S1)}] [Compute \,$s_\alpha$\, for \,$\alpha\in\Pi$\,.] 
\begin{itemize}
\item If \,$M$\, is of one of the types 
\,$\textsf{AIII}(n,q)$\, with \,$2q<n-1$\,, \,$\textsf{DIII}(n)$\, with \,$n$\, odd, and \,\textsf{EIII}\,, then put
$$ s_{\alpha_{\pi(1)}} := \tfrac{c_{\beta,\alpha_1}}{c_{\alpha_{\pi(1)},\beta}} \qmq{and} s_{\alpha} := 1 \text{ for all \,$\alpha \in \Pi \setminus \{\alpha_{\pi(1)}\}$\,} \; , $$
where \,$\alpha_k$\, and \,$\beta$\, have the same meaning as in Proposition~\ref{P:sigma:alphaadap}(a)(i).
\item If \,$M$\, is of any other type, then put \,$s_\alpha := 1$\, for all \,$\alpha\in\Pi$\,.
\end{itemize}
\item[\textbf{(S2)}] [Iterate on level.] Iterate steps (S3)--(S6) for \,$\ell=2,\dotsc,L$\,, where \,$L$\, denotes the maximal level of roots occurring in \,$\Delta$\,.
\item[\textbf{(S3)}] [Iterate on roots of level \,$\ell$\,.] Iterate steps (S4)--(S6) with \,$\alpha$\, running though all the roots in \,$\Delta$\, of level \,$\ell$\,.
\item[\textbf{(S4)}] [Find a decomposition of \,$\alpha$\,.] Let \,$\zeta,\eta$\, be positive roots so that \,$\alpha=\zeta+\eta$\, holds.
\item[\textbf{(S5)}] [Compute \,$s_{\alpha}$\,.] Put \,$s_\alpha := \tfrac{c_{\sigma\zeta,\sigma\eta}}{c_{\zeta,\eta}}\,s_\zeta\,s_\eta$\,.
\item[\textbf{(S6)}] (End of loops.)
\end{enumerate}

\begin{Remark}
For the decomposition \,$\alpha=\zeta+\eta$\, of the positive root \,$\alpha$\, requested in step (S4) of the algorithm, we can again use the decomposition
\,$\alpha=\zeta_\alpha+\eta_\alpha$\, already used in the computation of the Chevalley constants by algorithm (C), compare Remark~\ref{R:c:decomp}.
\end{Remark}

{\footnotesize
\emph{Proof for the correctness of the algorithm.}
Notice that each of the \,$s_\alpha$\, with \,$\alpha\in\Delta_+$\, is assigned to exactly once in the course of the algorithm, either in step (S1)
(if \,$\alpha$\, is simple), or in step (S5) (if \,$\alpha$\, is not simple). The assignments in step (S1) are correct by Proposition~\ref{P:sigma:alphaadap}(b)
(for the white roots) and Proposition~\ref{P:sigma:s}(d) (for the black roots). The correctness of the assignments in step (S5) follow by induction
on the level \,$\ell$\, via Proposition~\ref{P:sigma:s}(c); note that \,$\ell(\zeta),\ell(\eta)<\ell$\, holds in the situation of that step.
\strut\hfill $\Box$

}

\section{Formulas for the fundamental geometric tensors}
\label{Se:geom}

We are now ready to describe explicitly first the Lie bracket of \,$\lieg$\,, the inner product induced by its Killing form, and the involution \,$\sigma$\,,
thereby also the decomposition \,$\lieg = \liek \oplus \liem$\, induced by \,$\sigma$\,, and then the 
fundamental geometric tensors of \,$M$\,, namely the inner product on \,$\liem$\, and the curvature tensor \,$R$\, of \,$M$\,. 

We continue to use the notations of the preceding sections. In particular, we fix a Chevalley basis \,$(X_\alpha)$\, of \,$\lieg^{\C}$\, 
of the kind of Proposition~\ref{P:c:choice}, denote the corresponding Chevalley constants by \,$(c_{\alpha,\beta})$\, -- they can be calculated by algorithm (C) --
and consider the quantities \,$n_{kj}$\, and \,$s_\alpha$\, describing the 
involution \,$\sigma$\, as defined by Equations~\eqref{eq:sigma:sigmaalpha} and \eqref{eq:sigma:sdef}. We suppose that \,$\sigma$\, is of the kind
described in Proposition~\ref{P:sigma:alphaadap}(b), so that we have \,$s_{\alpha} \in \{\pm 1\}$\, for all \,$\alpha\in\Delta$\,, then the \,$s_\alpha$\,
can be calculated by algorithm (S). 

To describe the relevant tensors on \,$\lieg$\, (the Lie bracket, the inner product, and \,$\sigma$\,), it suffices to describe the behavior of \,$\liet$\,
and of the root spaces \,$\lieg_\alpha$\, with respect to these maps, because of the root space decomposition \eqref{eq:root:decomp-g}.

\begin{Prop}
\label{P:geom:g}
For any \,$\alpha\in\Delta$\, and \,$z \in \C$\, we put \,$V_\alpha(z) := \tfrac{1}{\sqrt{2}}\bigr(z\,X_\alpha - \overline{z}\,X_{-\alpha}\bigr)$\,. For formal reasons we put
\,$V_\alpha(z) := 0$\,, whenever \,$\alpha\in\liet^*$\, is a linear form which is not a root and \,$z\in \C$\,. 
Then we have for any \,$\alpha,\beta\in\Delta$\,:
\begin{enumerate}
\item \,$\lieg_\alpha = \Menge{V_\alpha(z)}{z\in\C}$\,.
\item For \,$H \in \liet$\, and \,$z,z' \in \C$\, we have
$$ [H,V_\alpha(z)] = V_\alpha(\alpha(H)\,z) $$
and
\begin{equation*}
[V_{\alpha}(z),V_\beta(z')] = \begin{cases}
\tfrac1{\sqrt{2}} \cdot \bigr(c_{\alpha,\beta} \, V_{\alpha+\beta}(z\,z') - c_{\alpha,-\beta} \,V_{\alpha-\beta}(z\,\overline{z'})\bigr) & \text{for \,$\beta \not\in\{\pm \alpha\}$\,} \\
\IM(\overline{z}\,z')\,i\alpha^\sharp & \text{for \,$\beta=\alpha$\,} \\
\IM(z\,z')\,i\alpha^\sharp & \text{for \,$\beta=-\alpha$\,} 
\end{cases} \; .
\end{equation*}
\item
\,$\liet$\, is orthogonal to \,$\lieg_\alpha$\, for every \,$\alpha\in \Delta$\,; \,$\lieg_\alpha$\, and \,$\lieg_\beta$\, are orthogonal for every \,$\alpha,\beta\in\Delta$\,
with \,$\alpha\not\in\{\pm \beta\}$\,; and besides the description of \,$\g{\cdot}{\cdot}$\, on \,$\liet\times\liet$\, in Section~\ref{Se:root} we have 
for every \,$\alpha\in\Delta$\,, \,$z,z'\in\C$\, 
$$ \g{V_\alpha(z)}{V_\alpha(z')} = \RE(z\cdot \overline{z'}) \; . $$
\item
Besides Equations~\eqref{eq:sigma:sigmaalpha} and \eqref{eq:sigma:sigmabeta}, 
which describe the action of \,$\sigma$\, on \,$\liet$\,, we have for any \,$\alpha\in\Delta$\, and \,$z\in\C$\,
$$ \sigma(V_\alpha(z)) = s_\alpha \cdot V_{\sigma(\alpha)}(z) \; . $$
\end{enumerate}
\end{Prop}

\beweis
\emph{For (a).} Because of the property \,$\overline{X_{\alpha}} = -X_{-\alpha}$\, of the Chevalley basis (see property (i) in Proposition~\ref{P:c:choice}),
we have 
\,$\Menge{V_\alpha(z)}{z\in \C} = (\lieg^{\C}_\alpha \oplus \lieg^{\C}_{-\alpha}) \cap \lieg = \lieg_\alpha$\,. 

\emph{For (b).} This is a straightforward computation, involving the definition of \,$V_\alpha(z)$\, and Equations~\eqref{eq:c:XaXb} and \eqref{eq:c:c-a-b}.

\emph{For (c).} The stated pairwise orthogonality of \,$\liet$\,, \,$\lieg_\alpha$\, and \,$\lieg_\beta$\, for \,$\beta \not\in \{\pm\alpha\}$\, is well-known,
and for \,$\alpha\in\Delta$\,, \,$z,z'\in\C$\, we have
$$ \g{V_\alpha(z)}{V_\alpha(z')} = -\tfrac12\,\vkap(z\,X_\alpha - \overline{z}\,X_{-\alpha} \,,\, z'\,X_\alpha - \overline{z'}\,X_{-\alpha})
\overset{(*)}{=} \tfrac12\,(z\,\overline{z'} + \overline{z}\,z') \, \vkap(X_\alpha,X_{-\alpha}) \overset{(\dagger)}{=} \RE(z\,\overline{z'}) \; , $$
where the equals sign marked \,$(*)$\, follows from the fact that \,$\vkap(X_\alpha,X_\alpha) = \vkap(X_{-\alpha},X_{-\alpha}) = 0$\, holds, and
the equals sign marked \,$(\dagger)$\, follows from Proposition~\ref{P:c:c}(b).

\emph{For (d).} This follows immediately from the definition of \,$V_{\alpha}(z)$\, and Equation~\eqref{eq:sigma:sdef}.
\beweisende

Let us now consider the decomposition \,$\lieg = \liek \oplus \liem$\, induced by \,$\sigma$\,. We have the root space decompositions
\begin{equation}
\label{eq:geom:decomp-km}
\liek = (\liet\cap\liek) \;\oplus \; \bigoplus_{\alpha\in \Delta_+^\sigma} \liek_\alpha \qmq{and}
\liem = \liea \;\oplus \; \bigoplus_{\alpha\in\Delta_+^\sigma} \liem_\alpha \;, 
\end{equation}
where \,$\Delta_+^\sigma \subset \Delta_+$\, is a subset such that for every \,$\alpha\in\Delta_+$\,, exactly one of the roots \,$\alpha,\sigma(\alpha),-\sigma(\alpha)$\,
(not necessarily pairwise unequal)
is a member of \,$\Delta_+^\sigma$\,, and where for any \,$\alpha\in(\liet^{\C})^*$\, we put \,$\liek_\alpha := (\lieg_\alpha + \lieg_{\sigma(\alpha)}) \cap \liek$\,
and \,$\liem_\alpha := (\lieg_\alpha + \lieg_{\sigma(\alpha)}) \cap \liem$\,.

To describe the fundamental geometric tensors of \,$M$\, on \,$\liem$\,, it again suffices to describe the behavior of \,$\liea$\, and the root spaces \,$\liem_\alpha$\,
with respect to them. We also describe  the behavior of the Lie bracket and the inner product with regard to \,$\liet\cap\liek$\, and the
root spaces \,$\liek_\alpha$\,. 

\begin{Prop}
\label{P:geom:M}
\begin{enumerate}
\item
Let \,$\alpha\in\Delta$\, be given.
\begin{enumerate}
\item
If \,$\sigma(\alpha) \not\in \{\pm \alpha\}$\, holds, we put for any \,$z\in\C$\,
$$ K_\alpha(z) := \tfrac{1}{\sqrt{2}}\,(V_\alpha(z)+s_\alpha\,V_{\sigma(\alpha)}(z)) \qmq{and}
M_\alpha(z) := \tfrac{1}{\sqrt{2}}\,(V_\alpha(z)-s_\alpha\,V_{\sigma(\alpha)}(z)) \; . $$
Then we have \,$\liek_\alpha = \Menge{K_\alpha(z)}{z\in\C}$\, and \,$\liem_\alpha = \Menge{M_\alpha(z)}{z\in\C}$\,.
\item
If \,$\sigma(\alpha)=\alpha$\, holds, we have
\,$\liek_\alpha = \Menge{V_\alpha(z)}{z\in\C}$\, and \,$\liem_\alpha = \{0\}$\,.
\item
If \,$\sigma(\alpha)=-\alpha$\, holds, we put for any \,$t\in\R$\,
$$ \wt{K}_\alpha(t) := \begin{cases}
V_\alpha(it) & \text{if \,$s_\alpha = 1$\,} \\
V_\alpha(t) & \text{if \,$s_\alpha = -1$\,} 
\end{cases} \qmq{and}
\wt{M}_\alpha(t) := \begin{cases}
V_\alpha(t) & \text{if \,$s_\alpha = 1$\,} \\
V_\alpha(it) & \text{if \,$s_\alpha = -1$\,} 
\end{cases} \; . $$
Then we have \,$\liek_\alpha = \Menge{\wt{K}_\alpha(t)}{t\in \R}$\, and \,$\liem_\alpha = \Menge{\wt{M}_\alpha(t)}{t\in \R}$\,.
\end{enumerate}
\item
For \,$\alpha,\beta\in\Delta$\, with \,$\beta\not\in\{\pm\alpha,\pm\sigma(\alpha)\}$\,, the spaces \,$(\liet\cap\liek), \liea, \liek_\alpha, \liem_\alpha,\liek_\beta,
\liem_\beta$\, are pairwise orthogonal. Moreover, for any \,$\alpha\in\Delta$\, we have
\begin{enumerate}
\item If \,$\sigma(\alpha) \not\in\{\pm \alpha\}$\,: For any \,$z,z'\in\C$\, we have \,$\g{K_\alpha(z)}{K_\alpha(z')} = \g{M_\alpha(z)}{M_\alpha(z')} 
= \RE(z\cdot \overline{z'})$\,.
\item If \,$\sigma(\alpha) = \alpha$\;: For any \,$z,z' \in \C$\, we have \,$\g{V_\alpha(z)}{V_\alpha(z')} = \RE(z\cdot \overline{z'})$\,.
\item If \,$\sigma(\alpha) = -\alpha$\;: For any \,$t,t' \in \R$\, we have \,$\g{\wt{K}_\alpha(t)}{\wt{K}_\alpha(t')} = \g{\wt{M}_\alpha(t)}{\wt{M}_\alpha(t')} 
= t\cdot t'$\,. 
\end{enumerate}
\item \,$\alpha,\beta\in\Delta$\,. To calculate the Lie bracket between elements of \,$\liek_\alpha\cup\liem_\alpha$\, and \,$\liek_\beta\cup\liem_\beta$\,,
we need to distinguish which of the three cases for \,$\alpha$\,, namely \,$\sigma(\alpha)=\alpha$\,, \,$\sigma(\alpha)=-\alpha$\,
or \,$\sigma(\alpha)\not\in\{\pm\alpha\}$\,, and similarly which case for \,$\beta$\, holds. By the combination of these cases, and use of the anti-symmetry
of the Lie bracket, there are in total six cases, which are handled separately in the following parts of the statement:
\begin{center}
\begin{tabular}{c|c|c|c}
& $\sigma(\alpha)\not\in\{\pm\alpha\}$ & \,$\sigma(\alpha)=\alpha$ & $\sigma(\alpha)=-\alpha$ \\
\hline
$\sigma(\beta)\not\in\{\pm\beta\}$ & (i) & (ii) & (iii) \\
$\sigma(\beta)=\beta$ & & (iv) & (v) \\
$\sigma(\beta)-\beta$ & & & (vi)
\end{tabular}
\end{center}
In the formulas, we have \,$z,z'\in\C$\, and \,$t\in\R$\,.
\begin{enumerate}
\item If \,$\sigma(\alpha) \not\in \{\pm \alpha\}$\, and \,$\sigma(\beta) \not\in \{\pm \beta\}$\, holds, we may suppose that either
\,$\beta \not\in \{\pm\alpha,\pm\sigma(\alpha)\}$\, or \,$\beta=\alpha$\, holds; indeed we can reduce the case \,$\beta \in \{\pm\alpha,\pm\sigma(\alpha)\}$\,
to \,$\beta=\alpha$\, by application of the formulas
$$ K_{\sigma(\beta)}(z') = s_\beta\,K_\beta(z'),\; M_{\sigma(\beta)}(z') = -s_\beta\,M_\beta(z'), \qmq{and}
K_{-\beta}(z') = -K_\beta(\overline{z'}),\; M_{-\beta}(z') = -M_\beta(\overline{z'}) \; . $$
Then we have:
{\tiny
\begin{align*}
[K_\alpha(z),K_\beta(z')] & = \begin{cases}
\tfrac12 \cdot \bigr( c_{\alpha,\beta}\,K_{\alpha+\beta}(z\,z') - c_{\alpha,-\beta}\,K_{\alpha-\beta}(z\,\overline{z'})
        + s_\alpha\,c_{\sigma(\alpha),\beta}\,K_{\sigma(\alpha)+\beta}(z\,z') - s_\alpha\,c_{\sigma(\alpha),-\beta}\,K_{\sigma(\alpha)-\beta}(z\,\overline{z'})\bigr)
& \text{for \,$\beta\not\in\{\pm\alpha,\pm\sigma(\alpha)\}$\,} \\
\tfrac12\,\IM(\overline{z}\,z')\,i\,(\alpha^\sharp + \sigma(\alpha)^\sharp) + \sqrt{2}\,s_\alpha\,c_{\alpha,-\sigma(\alpha)}\,\wt{K}_{\alpha-\sigma(\alpha)}(\IM(\overline{z}\,z'))
        & \text{for \,$\beta=\alpha$\,}
\end{cases} \\
[K_\alpha(z),M_\beta(z')] & = \begin{cases}
\tfrac12 \cdot \bigr( c_{\alpha,\beta}\,M_{\alpha+\beta}(z\,z') - c_{\alpha,-\beta}\,M_{\alpha-\beta}(z\,\overline{z'})
        + s_\alpha\,c_{\sigma(\alpha),\beta}\,M_{\sigma(\alpha)+\beta}(z\,z') - s_\alpha\,c_{\sigma(\alpha),-\beta}\,M_{\sigma(\alpha)-\beta}(z\,\overline{z'})\bigr)
& \text{for \,$\beta\not\in\{\pm\alpha,\pm\sigma(\alpha)\}$\,} \\
\tfrac12\,\IM(\overline{z}\,z')\,i\,(\alpha^\sharp - \sigma(\alpha)^\sharp) + \sqrt{2}\,s_\alpha\,c_{\alpha,-\sigma(\alpha)}\,\wt{M}_{\alpha-\sigma(\alpha)}(\RE(\overline{z}\,z'))
        & \text{for \,$\beta=\alpha$\,}
\end{cases}\\
[M_\alpha(z),M_\beta(z')] & = \begin{cases}
\tfrac12 \cdot \bigr( c_{\alpha,\beta}\,K_{\alpha+\beta}(z\,z') - c_{\alpha,-\beta}\,K_{\alpha-\beta}(z\,\overline{z'})
        - s_\alpha\,c_{\sigma(\alpha),\beta}\,K_{\sigma(\alpha)+\beta}(z\,z') + s_\alpha\,c_{\sigma(\alpha),-\beta}\,K_{\sigma(\alpha)-\beta}(z\,\overline{z'})\bigr)
& \text{for \,$\beta\not\in\{\pm\alpha,\pm\sigma(\alpha)\}$\,} \\
\tfrac12\,\IM(\overline{z}\,z')\,i\,(\alpha^\sharp + \sigma(\alpha)^\sharp) - \sqrt{2}\,s_\alpha\,c_{\alpha,-\sigma(\alpha)}\,\wt{K}_{\alpha-\sigma(\alpha)}(\IM(\overline{z}\,z'))
        & \text{for \,$\beta=\alpha$\,}
\end{cases}
\end{align*}
}
\item
If \,$\sigma(\alpha)=\alpha$\, and \,$\sigma(\beta)\not\in\{\pm\beta\}$\, holds, we have
\begin{align*}
[V_\alpha(z),K_\beta(z')] & = c_{\alpha,\beta}\,K_{\alpha+\beta}(z\,z') - c_{\alpha,-\beta}\,K_{\alpha-\beta}(z\,\overline{z'}) \\
[V_\alpha(z),M_\beta(z')] & = c_{\alpha,\beta}\,M_{\alpha+\beta}(z\,z') - c_{\alpha,-\beta}\,M_{\alpha-\beta}(z\,\overline{z'}) \; . 
\end{align*}
\item
If \,$\sigma(\alpha)=-\alpha$\, and \,$\sigma(\beta)\not\in\{\pm\beta\}$\, holds, we put \,$\zeta := i$\, if \,$s_\alpha=1$\,, \,$\zeta := 1$\, if \,$s_\alpha=-1$\,,
then we have
\begin{align*}
[\wt{K}_\alpha(t),K_\beta(z)] & = c_{\alpha,\beta}\,K_{\alpha+\beta}(\zeta t\,z) - c_{\alpha,-\beta}\,K_{\alpha-\beta}(\zeta t\,\overline{z}) \\
[\wt{K}_\alpha(t),M_\beta(z)] & = c_{\alpha,\beta}\,M_{\alpha+\beta}(\zeta t\,z) - c_{\alpha,-\beta}\,M_{\alpha-\beta}(\zeta t\,\overline{z}) \\
[\wt{M}_\alpha(t),K_\beta(z)] & = c_{\alpha,\beta}\,M_{\alpha+\beta}((1+i-\zeta)t\,z) - c_{\alpha,-\beta}\,M_{\alpha-\beta}((1+i-\zeta)t\,\overline{z}) \\
[\wt{M}_\alpha(t),M_\beta(z)] & = c_{\alpha,\beta}\,K_{\alpha+\beta}((1+i-\zeta)t\,z) - c_{\alpha,-\beta}\,K_{\alpha-\beta}((1+i-\zeta)t\,\overline{z}) \; . 
\end{align*}
\item
If \,$\sigma(\alpha)=\alpha$\, and \,$\sigma(\beta)=\beta$\, holds, then \,$[V_{\alpha}(z),V_\beta(z')]$\, is given by Proposition~\ref{P:geom:g}(b).
\item
If \,$\sigma(\alpha)=-\alpha$\, and \,$\sigma(\beta)=\beta$\, holds, we again put \,$\zeta := i$\, if \,$s_\alpha=1$\,, \,$\zeta := 1$\, if \,$s_\alpha=-1$\,,
then we have
\begin{align*}
[\wt{K}_\alpha(t),V_\beta(z)] & = \sqrt{2}\,c_{\alpha,\beta}\,K_{\alpha+\beta}(\zeta t\,z) \\
[\wt{M}_\alpha(t),V_\beta(z)] & = \sqrt{2}\,c_{\alpha,\beta}\,M_{\alpha+\beta}(\zeta t\,z) \; . 
\end{align*}
\item
If \,$\sigma(\alpha)=-\alpha$\, and \,$\sigma(\beta)=-\beta$\, holds, then the values of \,$[\wt{K}_\alpha(t),\wt{K}_\beta(s)], [\wt{K}_\alpha(t),\wt{M}_\beta(s)]$\,
and \,$[\wt{M}_\alpha(t),\wt{M}_\beta(s)]$\, are given by Proposition~\ref{P:geom:g}(b) in conjunction with
the definition of \,$\wt{K}_\alpha(t)$\,, \,$\wt{M}_\alpha(t)$\,.
\end{enumerate}
\item
The curvature tensor \,$R$\, of \,$M$\, is given by
$$ R(u,v)w = -[[u,v],w] \qmq{for \,$u,v,w\in\liem$\,} $$
and can therefore be calculated by two-fold application of (c).
\end{enumerate}
\end{Prop}

\beweis
\emph{For (a).} In the case of (i) we have \,$\sigma(V_\alpha(z)) = s_\alpha\,V_{\sigma(\alpha)}(z)$\,, in the case of (ii) we have 
\,$s_\alpha=1$\, by Proposition~\ref{P:sigma:s}(d) and therefore \,$\sigma(V_\alpha(z))=V_\alpha(z)$\,, and in the case of (iii) we have
\,$\sigma(V_\alpha(z)) = s_\alpha\,V_{-\alpha}(z) = -s_\alpha\,V_\alpha(\overline{z})$\,. Therefrom the statements follow.

\emph{For (b).} These statements follow by application of Proposition~\ref{P:geom:g}(c) with the explicit descriptions of \,$\liek_\alpha$\, and \,$\liem_\alpha$\,
from (a).

\emph{For (c).} The formulas for the Lie bracket are derived from the definitions in (a), the equation in Proposition~\ref{P:geom:g}(b) for \,$[V_\alpha(z),V_\beta(z')]$\,
and the computational rules for the \,$c_{\alpha,\beta}$\, (Proposition~\ref{P:c:c}) and the \,$s_\alpha$\, (Proposition~\ref{P:sigma:s}). For (i) one also needs
to use the fact that for any \,$\alpha\in\Delta$\, we have \,$\alpha+\sigma(\alpha)\not\in\Delta$\, (see \cite{Loos:1969-2}, Proposition~VI.3.3(c), p.~73).

\emph{(d)} is a well-known statement.
\beweisende

\section{An application: Totally geodesic submanifolds in \,$\SU(3)/\SO(3)$\,}
\label{Se:AI(2)}

As an application of the preceding construction of the curvature tensor, we show how to classify the Lie triple systems (i.e.~those linear subspaces of \,$\liem$\, which are
invariant under the curvature tensor) in the Riemannian symmetric space \,$\SU(3)/\SO(3)$\,. They are exactly the tangent spaces of totally geodesic submanifolds
of \,$\SU(3)/\SO(3)$\, passing through the origin point.

Note that the rank of the symmetric space \,$M:=\SU(3)/\SO(3)$\, is \,$2$\,, and thus equals the rank of its transvection group \,$G:=\SU(3)$\,. If we consider
the splitting \,$\lieg := \liesu(3) = \liek \oplus\liem$\, induced by the symmetric structure of \,$M$\,, any maximal flat subspace \,$\liea$\, of \,$\liem$\,
therefore already is a Cartan subalgebra of \,$\lieg$\,, and thus the root system \,$\Delta$\, of \,$\lieg$\, equals the ``restricted'' root system
of \,$\liem$\,. In the present case, this root system is of type \,$A_2$\,, so if \,$\{\alpha_1,\alpha_2\}$\, is a system of simple roots in \,$\Delta$\,, we have
\,$\Delta = \{\pm \alpha_1,\pm \alpha_2,\pm \alpha_3\}$\, with \,$\alpha_3 := \alpha_1+\alpha_2$\,. Moreover, we have \,$\sigma(\alpha_k) = -\alpha_k$\, and
thus the root space decomposition
\begin{equation}
\label{eq:AI(2):m-decomp}
\liem = \liea \oplus \bigoplus_{k=1}^3 \liem_{\alpha_k} \qmq{with} \liem_{\alpha_k} = \R\,M_{\alpha_k}(1) \qmq{for \,$k=1,2,3$\,} 
\end{equation}
(see Equation~\eqref{eq:geom:decomp-km} and Proposition~\ref{P:geom:M}(a)(iii)). For every \,$k\in\{1,2,3\}$\,, the linear form \,$\alpha_k$\, is purely imaginary
on \,$\liea$\, because the symmetric space \,$M$\, is of compact type, and thus we have \,$H_k := (\tfrac1i\,\alpha_k)^\sharp \in \liea$\,. 

\enlargethispage{2em} 

\begin{Prop}
\label{P:AI(2):AI(2)}
Let \,$\{0\} \neq \liem' \subsetneq \liem$\, be a linear subspace. Then \,$\liem'$\, is a Lie triple system in \,$\liem$\, if and only if there exists
a maximal flat subspace \,$\liea \subset \liem$\, and a system of simple roots \,$\{\alpha_1,\alpha_2\}$\, in the root system \,$\Delta$\, of \,$\liem$\,
(or of \,$\lieg$\,) with respect to \,$\liea$\,, so that with \,$\alpha_3 := \alpha_1+\alpha_2$\, and \,$H_k$\, defined in relation to these \,$\alpha_k$\,
as above, \,$\liem'$\, is of one of the following types:
\begin{itemize}
\item[(G)] \,$\liem' = \R\,H$\, with some \,$H \in \liea$\,
\item[(T)] \,$\liem' = \liea$\,
\item[(S)] \,$\liem' = \R\,H_3 \oplus \R\,M_{\alpha_3}(1)$\,
\item[(M)] \,$\liem' = \R\,H_3 \oplus \R\,(M_{\alpha_1}(1) + M_{\alpha_2}(1))$\,
\item[(P)] \,$\liem' = \liea \oplus \R\,M_{\alpha_3}(1)$\,
\end{itemize}
Two Lie triple systems of types other than (G) are congruent to each other under the isotropy action on \,$\liem$\, if and only if they are of the same type. 
A Lie triple system is maximal if and only if it is either of type (P) or of type (M).
\end{Prop}

\begin{Remark}
The totally geodesic submanifolds corresponding to Lie triple systems of type (M) --- they are isometric to an \,$\RP^2$\, of sectional curvature
\,$\tfrac12$\, --- are missing from the 
classification of maximal totally geodesic submanifolds in Riemannian symmetric spaces of rank \,$2$\, by \textsc{Chen} and \textsc{Nagano} in 
\cite{Chen/Nagano:totges2-1978}. (Notice that they are not contained in the local sphere products corresponding to type (P).)
\end{Remark}

{\footnotesize\emph{Proof of Proposition~\ref{P:AI(2):AI(2)}.}
First we note that the spaces described in the proposition are indeed Lie triple systems: For the types (G) and (T) this is true because they are flat,
for the types (P) and (S) it is true because they correspond to the closed root subsystem \,$\{\pm \alpha_3\}$\, of \,$\Delta$\,, for type (M): In
the setting of that type we put
\,$H := H_3 \in \liea$\, and \,$v := M_{\alpha_1}(1) + M_{\alpha_2}(1)$\,, then we have \,$\alpha_1(H) = \alpha_2(H)$\,, hence \,$R(H,v)H \in \R\,v \subset \liem'$\,;
we also see \,$R(H,v)v \in \R\,H \subset \liem'$\, by explicitly calculating this vector via the description of \,$R$\, developed 
in this paper (see Equation~\eqref{eq:AI(2):AI(2):RHvv} below), and therefore \,$\liem' = \R\,H\oplus\R\,v$\, is a Lie triple system. Presuming that the
list of Lie triple systems given in the proposition is complete, we see that the types (P) and (M) are maximal, whereas the inclusions
(G) $\subset$ (T)  $\subset$ (P) and (S)  $\subset$ (P) hold, showing that no other types are maximal.

It remains to show that any given Lie triple system \,$\liem'$\, of \,$\liem$\, is of one of the types given in the proposition. 
We will base the proof on the fact that also for \,$\liem'$\, we have a root space decomposition, and on the relations that hold between
that decomposition and the root space decomposition \eqref{eq:AI(2):m-decomp} for \,$\liem$\,; for a detailed description of these relations, see
\cite{Klein:2007-claQ}, Section~2. We have \,$\rk(\liem') \leq \rk(\liem) = 2$\,,
and therefore \,$\rk(\liem') \in \{1,2\}$\,. We will handle the two possibilities for the rank of \,$\liem'$\, separately. 

If \,$\rk(\liem')=2$\,, then any maximal flat subspace \,$\liea$\, of \,$\liem'$\, also is a maximal flat subspace of \,$\liem$\,,
and if we denote the root systems of \,$\liem'$\, resp.~of \,$\liem$\, with respect to \,$\liea$\, by \,$\Delta'$\, resp.~\,$\Delta$\,, then \,$\Delta' 
\subset \Delta$\, holds; moreover if we consider the root space decomposition \,$\liem' = \liea \oplus \bigoplus_{\alpha\in\Delta'_+}\liem_{\alpha}'$\, of \,$\liem'$\,,
then we have \,$\{0\} \neq \liem_\alpha' \subset \liem_\alpha$\, and therefore
\,$\liem_\alpha' = \liem_\alpha$\, for any \,$\alpha\in\Delta'$\,, because \,$\liem_\alpha$\, is 1-dimensional.

There are thus three possibilities: \,$\Delta' = \varnothing$\,, \,$\Delta' = \{\pm \alpha_1\}$\, for some \,$\alpha_1 \in \Delta$\,,
and \,$\Delta' = \Delta$\,. If \,$\Delta' = \varnothing$\, holds, we have \,$\liem' = \liea$\,, and thus \,$\liem'$\, is of type (T). If \,$\Delta' = \{\pm\alpha_1\}$\,
holds, then we have \,$\liem' = \liea\oplus\liem_{\lambda_1}$\,, and thus \,$\liem'$\, is of type (P). Finally, \,$\Delta'=\Delta$\, is possible only for \,$\liem'=\liem$\,. 

If \,$\rk(\liem')=1$\,, we fix some \,$H\in\liem'\setminus\{0\}$\,, then \,$\liea' := \R\,H$\, is a maximal flat subspace of \,$\liem'$\,. We choose a maximal
flat subspace \,$\liea$\, of \,$\liem$\, with \,$\liea' = \liea\cap\liem'$\,. If \,$\dim(\liem')=1$\, holds, we have \,$\liem' = \liea'$\,, and therefore
\,$\liem'$\, then is of type (G). Thus we suppose \,$\dim(\liem') \geq 2$\, in the sequel.
Then we have (see \cite{Klein:2007-claQ}, Proposition~2.3) either \,$H \in \R\,(\tfrac1i\,\alpha)^\sharp$\, for some \,$\alpha\in\Delta$\,,
or \,$H \perp \tfrac1i(\alpha^\sharp - \beta^\sharp)$\, for some \,$\alpha,\beta\in\Delta$\,, \,$\alpha\neq\beta$\,. It follows that there exists a system of simple
roots \,$\{\alpha_1,\alpha_2\}$\, of \,$\Delta$\, so that with \,$\alpha_3 := \alpha_1+\alpha_2\in\Delta$\, we have (after scaling \,$H$\, appropriately) 
either \,$H = H_3$\, or \,$H = H_3 + \tfrac12\,(H_1-H_3) = H_1 + \tfrac12\,H_2$\,. 
We will treat these two possibilities separately below, but in either case we have
\,$\varnothing \neq \Delta' \subset \Mengegr{\alpha|\liea'}{\alpha\in\Delta,\, \alpha|\liea'\neq 0}$\, and the following root space decomposition for \,$\liem'$\,:
\begin{equation}
\label{eq:AI(2):AI(2):m'-decomp}
\liem' = \liea' \oplus \bigoplus_{\lambda\in\Delta'_+} \liem_\lambda' \qmq{with} \liem_\lambda' = \left( \bigoplus_{\substack{\alpha \in \Delta \\ \alpha|\liea' = \lambda}} 
\liem_\alpha \right) \cap \liem' \qmq{for \,$\lambda\in\Delta'$\,} \; .
\end{equation}

Let us now first consider the case \,$H = H_3$\,. Then we have \,$\Delta' \subset \{\pm \lambda,\pm 2\lambda\}$\, with \,$\lambda := \alpha_1|\liea'
= \alpha_2|\liea'$\,; note \,$2\lambda = \alpha_3|\liea'$\,. By Equation~\eqref{eq:AI(2):AI(2):m'-decomp} we therefore have 
\begin{equation}
\label{eq:AI(2):AI(2):subspaces}
\liem_\lambda' \subset \liem_{\alpha_1} \oplus \liem_{\alpha_2} \qmq{and} \liem_{2\lambda}' \subset \liem_{\alpha_3} \; . 
\end{equation}
Let \,$v \in \liem_\lambda'$\, be given, say \,$v = M_{\alpha_1}(a) + M_{\alpha_2}(b)$\, with \,$a,b\in\R$\, (see Equation~\eqref{eq:AI(2):AI(2):subspaces}).
Because \,$\liem'$\, is a Lie triple system, we have \,$R(H,v)v \in \liem'$\,, and using the representation
of the curvature tensor \,$R$\, given in the present paper, we can actually calculate that vector explicitly (most easily by 
using the \textsf{Maple} implementation of the algorithms mentioned in the Introduction): 
\begin{equation}
\label{eq:AI(2):AI(2):RHvv}
R(H,v)v = \tfrac12\,(a^2\,H_1 + b^2\, H_2) \; .
\end{equation}
Therefore \,$R(H,v)v$\, is a member of \,$\liem' \cap \liea = \liea' = \R\,H = \R\,H_3 = \R\,(H_1+H_2)$\,, and hence
\eqref{eq:AI(2):AI(2):RHvv} gives \,$a^2 = b^2$\, and thus \,$a=\pm b$\,. It follows that we have either \,$\liem_\lambda' \subset \R\,(M_{\alpha_1}(1) + M_{\alpha_2}(1))$\,
or \,$\liem_\lambda' \subset \R\,(M_{\alpha_1}(1) - M_{\alpha_2}(1))$\,. In fact, we can suppose without loss of generality
\begin{equation}
\label{eq:AI(2):AI(2):mlambda'}
\liem_\lambda' \subset \R\,(M_{\alpha_1}(1) + M_{\alpha_2}(1)) \; . 
\end{equation}

Let us now suppose \,$2\lambda\in\Delta'$\,, then we have \,$\liem_{2\lambda}' = \liem_{\alpha_3}$\, by \eqref{eq:AI(2):AI(2):subspaces}, and therefore
\,$M_{\alpha_3}(1) \in \liem'$\, holds. Again using our representation of \,$R$\,, we calculate the element \,$R(M_{\alpha_3}(1),v)v$\, of \,$\liem'$\,:
\begin{equation}
\label{eq:AI(2):AI(2):RM3vv}
R(M_{\alpha_3}(1),v)v = \tfrac12 ab\,(H_1-H_2) + \tfrac14\,M_{\alpha_3}(a^2 + b^2) \; .
\end{equation}
Because the \,$\liea$-component of this vector, \,$\tfrac12 ab\,(H_1-H_2)$\,, is contained in \,$\liea' = \R\,(H_1+H_2)$\,,
we see that \,$ab=0$\, holds. Because we have already seen \,$b=\pm a$\,, \,$a=b=0$\, and thus \,$v=0$\, follows. So we have shown that 
\,$2\lambda\in\Delta$\, implies \,$\lambda\not\in\Delta$\,. Thus Equations~\eqref{eq:AI(2):AI(2):m'-decomp} and \eqref{eq:AI(2):AI(2):subspaces} imply
that in this case \,$\liem' = \R\,H_3 \oplus \liem_{\alpha_3}$\, holds, and therefore \,$\liem'$\, is of type (S).

On the other hand, for \,$2\lambda\not\in\Delta'$\, we have \,$\lambda \in \Delta'$\, and thus by Equation~\eqref{eq:AI(2):AI(2):mlambda'}:
\,$\liem_\lambda' = \R\,(M_{\alpha_1}(1) + M_{\alpha_2}(1))$\,. It now follows from Equation~\eqref{eq:AI(2):AI(2):m'-decomp} that \,$\liem'$\,
is of type (M).

Finally, we consider the case \,$H = H_1 + \tfrac12\,H_2$\,. Then we have \,$\alpha_1(H) = \alpha_3(H) > 0$\, and \,$\alpha_2(H)=0$\,,
and thus \,$\Delta' \subset \{\pm \lambda\}$\, with \,$\lambda := \alpha_1|\liea' = \alpha_3|\liea'$\, and \,$\liem_\lambda' \subset \liem_{\alpha_1} \oplus \liem_{\alpha_3}$\,.
Let \,$v \in \liem_\lambda'$\, be given, say \,$v = M_{\alpha_1}(a) + M_{\alpha_3}(b)$\, with \,$a,b \in \R$\,. We once again use our representation of \,$R$\,
to compute the element \,$R(H,v)v$\, of \,$\liem'$\, explicitly:
\begin{equation}
\label{eq:AI(2):AI(2):RHvv-2}
R(H,v)v = \tfrac34\,\bigr(\,(a^2+b^2)\,H_1 + b^2\,H_2\,\bigr) - \tfrac34\,M_{\alpha_2}(ab) \; . 
\end{equation}
Because \,$\liem'$\, is orthogonal to \,$\liem_{\alpha_2}$\, by Equation~\eqref{eq:AI(2):AI(2):m'-decomp}, 
the \,$\liem_{\alpha_2}$-component of \eqref{eq:AI(2):AI(2):RHvv-2} must be zero, so we have \,$ab=0$\,.
Also, the \,$\liea$-component of \eqref{eq:AI(2):AI(2):RHvv-2} must be a member of \,$\liea' = \R\,(H_1 + \tfrac12\,H_2)$\,, so we have
\,$a^2 + b^2 = 2b^2$\,, hence \,$a = \pm b$\,. From these two equations, \,$a=b=0$\, and hence \,$v=0$\, follows. So we have \,$\liem_\lambda' = \{0\}$\,
and thus \,$\liem' = \liea'$\, is 1-dimensional. This shows that the case \,$H = H_1 + \tfrac12\,H_2$\, in fact cannot occur
(for \,$\dim(\liem') \geq 2$\,), and this completes the classification.
\strut\hfill$\Box$

}

We finally discuss the totally geodesic submanifolds corresponding to the various types of Lie triple systems of \,$\SU(3)/\SO(3)$\, found in 
Proposition~\ref{P:AI(2):AI(2)}. 
For the purpose of describing the metric of the submanifolds, we suppose that the \,$\SU(3)$-invariant Riemannian
metric on \,$\SU(3)/\SO(3)$\, is the one induced by the usual inner product on \,$\End(\C^3)$\,, namely the one given by \,$\g{A}{B} := \RE(\tr(B^* A))$\, for
\,$A,B\in\End(\C^3)$\,. Then the root vectors \,$H_k$\, have length \,$\sqrt{2}$\,.
With this choice of Riemannian metric on \,$\SU(3)/\SO(3)$\,,
the totally geodesic submanifolds corresponding to the various types of Lie triple systems have the following isometry type:

{\scriptsize
\begin{center}
\begin{tabular}{|c||c|c|c|c|c|}
\hline
type of Lie triple system & (G) & (T) & (S) & (M) & (P) \\
\hline
isometry type & \,$\R$\, or \,$S^1_r$\, & \,$(S^1_{r=\sqrt{3/2}} \times S^1_{r=\sqrt{1/2}})/\{\pm \id\}$\, & \,$S^2_{r=\sqrt{1/2}}$\, & \,$\RP^2_{\vkap=1/2}$\, & \,$(S^1_{r=\sqrt{3/2}} \times S^2_{r=\sqrt{1/2}})/\{\pm \id\}$\, \\
\hline
\end{tabular}
\end{center}

}

\medskip

The totally geodesic submanifolds of type (G) are the traces of geodesics in \,$\SU(3)/\SO(3)$\,, and the submanifolds of type (T) are the
maximal flat tori in \,$\SU(3)/\SO(3)$\,. A totally geodesic submanifold of \,$\SU(3)/\SO(3)$\, is reflective (i.e.~is a connected component of the 
fixed point set of an involutive isometry, see for example \cite{Leung:reflective-1979} and other papers by \textsc{Leung}) 
if and only if it is either of type (M) or of type (P). For a Lie triple system corresponding to a reflective
submanifold, the orthogonal complement is again a Lie triple system, and in this way, the types (M) and (P) correspond to each other. In fact, the
submanifolds of type (M) are polars (i.e.~they are connected components of the fixed point set of the geodesic symmetry of \,$\SU(3)/\SO(3)$\,, see 
\cite{Chen/Nagano:totges2-1978}, \S 2) and the submanifolds of type (P) are the corresponding meridians.

Note that via these concepts, all totally geodesic submanifolds of \,$\SU(3)/\SO(3)$\, can be obtained in a ``natural'' way: The submanifolds of type (M)
are the polars of \,$\SU(3)/\SO(3)$\,, and the submanifolds of type (P) correspond to them as meridians. The remaining totally geodesic submanifolds,
namely those of type (S), (T) and (G), are obtained as the ``obvious'' totally geodesic submanifolds of the meridians.

To prove that the totally geodesic submanifolds are indeed of the isometry types given in the above table, and to describe their position 
in \,$\SU(3)/\SO(3)$\,, we now give totally geodesic embeddings for each type of totally geodesic submanifold explicitly.

\textbf{Type (S).}
Consider the Lie group monomorphism
$$ \Phi_0: \SU(2) \to \SU(3), \; B \mapsto \left( \begin{matrix} B & 0 \\ 0 & 1 \end{matrix} \right) \; , $$
which is isometric with regard to the Riemannian metrics induced on \,$\SU(2)$\, resp.~\,$\SU(3)$\, by the usual endomorphism inner product.
We have \,$\Phi_0^{-1}(\SO(3)) = \SO(2)$\,, and \,$\Phi_0$\, is compatible with the Lie group involutions induced by the symmetric space structures of
\,$\SU(2)/\SO(2)$\, resp.~\,$\SU(3)/\SO(3)$\,. Therefore \,$\Phi_0$\, gives rise to a totally geodesic isometric embedding 
$$ \underline{\Phi}_0: \SU(2)/\SO(2) \to \SU(3)/\SO(3) \; , $$
its image is a totally geodesic submanifold of \,$\SU(3)/\SO(3)$\,, which turns out to be of type (S). 
Thus the totally geodesic submanifolds of type (S) are isometric to \,$\SU(2)/\SO(2)$\,. 

\,$\SU(2)/\SO(2)$\, is a simply connected, 2-dimensional, irreducible Riemannian symmetric space of compact type, and hence isometric 
to a 2-sphere of some specific radius \,$r$\,. The curve \,$\gamma: \R \to \SU(2)/\SO(2),\;
t \mapsto \left( \begin{smallmatrix} e^{it/\sqrt{2}} & 0 \\ 0 & e^{-it/\sqrt{2}} \end{smallmatrix}\right)\cdot \SO(2)$\,
is a unit speed geodesic of \,$\SU(2)/\SO(2)$\, with period \,$\sqrt{2}\,\pi = 2\pi r$\,. Therefore \,$\SU(2)/\SO(2)$\, 
(and hence, any totally geodesic submanifold of \,$\SU(3)/\SO(3)$\, of type (S)) is isometric to
\,$S^2_{r=1/\sqrt{2}}$\,. 

\textbf{Type (P).}
To describe a totally geodesic embedding of type (P), we ``extend'' the embedding \,$\Phi_0$\, described above in the following way:
$$ \Phi: S^1_{r=\sqrt{6}} \times \SU(2) \to \SU(3),  \; (\lambda,B) \mapsto \begin{pmatrix} \tfrac{\lambda}{\sqrt{6}}\,B & 0 \\ 0 & \tfrac{6}{\lambda^{2}} \end{pmatrix} \; , $$
where we regard \,$S^1_{r=\sqrt{6}}$\, as the circle \,$\Mengegr{z\in\C}{|z|^2=6}$\, in \,$\C$\,. Note that \,$\Phi(\sqrt{6},\,\cdot\,) = \Phi_0$\, holds.
The differential of \,$\Phi$\, at \,$(\sqrt6, \id) \in S^1_{r=\sqrt{6}} \times \SU(2)$\, is given by
$$ T_{\sqrt{6}} S^1_{r=\sqrt{6}} \times \liesu(2) \to \liesu(3),\; 
(it, X) \mapsto \begin{pmatrix} \tfrac{it}{\sqrt{6}}\,\id + X & 0 \\ 0 & - \tfrac{2\,it}{\sqrt{6}} \end{pmatrix} \;, $$
where we identify the tangent space of \,$S^1_{r=\sqrt{6}}$\, at \,$\sqrt{6}$\, with \,$i\R$\,. Using this presentation of the differential of \,$\Phi$\,
it is easy to see that \,$\Phi$\, is isometric.
Moreover, we have \,$\Phi^{-1}(\SO(3)) = K \cup gK$\, with \,$K := \{\pm \sqrt{6}\} \times \SO(2)$\, 
and \,$g := (\sqrt{6}\,i,J)$\, where \,$J := \left( \begin{smallmatrix} i & 0 \\ 0 & -i \end{smallmatrix} \right) \in \SU(2)$\,, 
and again \,$\Phi$\, is compatible with the
symmetric structures of the symmetric spaces involved. Therefore \,$\Phi$\, gives rise to a totally geodesic, isometric embedding
$$ \underline{\Phi}: (S^1_{r=\sqrt{6}} \times \SU(2)) / (K\cup gK) \to \SU(3)/\SO(3) \; . $$
The image of \,$\underline{\Phi}$\, is a totally geodesic submanifold of \,$\SU(3)/\SO(3)$\,, which turns out to be of type (P). Therefore the totally geodesic
submanifolds of \,$\SU(3)/\SO(3)$\, of type (P) are isometric to \,$(S^1_{r=\sqrt{6}} \times \SU(2)) / (K \cup gK)$\,. 

It remains to describe the isometry type of \,$(S^1_{r=\sqrt{6}} \times \SU(2)) / (K \cup gK)$\, more succinctly. 
We have \,$(S^1_{r=\sqrt{6}} \times \SU(2)) / K = (S^1_{r=\sqrt{6}}/\{\pm \sqrt{6}\}) \times (\SU(2)/\SO(2))$\,, where \,$S^1_{r=\sqrt{6}}/\{\pm \sqrt{6}\}$\,
is isometric to \,$S^1_{r=\sqrt{6}/2}$\,, and \,$\SU(2)/\SO(2)$\, is isometric to \,$S^2_{r=1/\sqrt{2}}$\, as we saw in the treatment of type (S).
Hence \,$(S^1_{r=\sqrt{6}} \times \SU(2)) / K$\, is isometric to \,$S^1_{r=\sqrt{6}/2} \times S^2_{r=1/\sqrt{2}}$\,. Because \,$K,gK \in (S^1_{r=\sqrt{6}} \times \SU(2)) / K$\,
correspond to a pair of antipodal points in \,$S^1_{r=\sqrt{6}/2} \times S^2_{r=1/\sqrt{2}}$\, under this isometry, it follows that 
\,$(S^1_{r=\sqrt{6}} \times \SU(2)) / (K \cup gK)$\, (and hence, any totally geodesic submanifold of \,$\SU(3)/\SO(3)$\, of type (P)) 
is isometric to \,$(S^1_{r=\sqrt{6}/2} \times S^2_{r=1/\sqrt{2}})/\{\pm \id\}$\,. 

\textbf{Type (T).}
The Lie triple systems of type (T) are the maximal flat subspaces of \,$\liem$\,, so the corresponding totally geodesic submanifolds are 
the maximal flat tori of \,$\SU(3)/\SO(3)$\,. These Lie triple systems are contained in Lie triple systems of type (P)
(as can be seen by the explicit description of the Lie triple systems in Proposition~\ref{P:AI(2):AI(2)}), and therefore
the maximal flat tori of \,$\SU(3)/\SO(3)$\, are contained in totally geodesic submanifolds of type (P). A totally geodesic, isometric embedding
of type (T) can be obtained by fixing a one-parameter subgroup \,$C$\, of \,$\SU(2)$\, which runs orthogonal to \,$\SO(2)\subset\SU(2)$\, 
and restricting \,$\underline{\Phi}$\, to \,$(S^1_{r=\sqrt{6}} \times C) / (K\cup gK)$\,. Therefore the maximal flat tori of \,$\SU(3)/\SO(3)$\, 
are isometric to \,$(S^1_{r=\sqrt{6}/2} \times S^1_{r=1/\sqrt{2}})/\{\pm \id\}$\,. 

To understand the geometry of the maximal tori in \,$\SU(3)/\SO(3)$\, better, 
we consider the lattice \,$\wt{\Gamma} := \Z\,(\sqrt{6}\,\pi,0) \oplus \Z\,(0,\sqrt{2}\,\pi)$\, in \,$\R^2$\,. Then
\,$S^1_{r=\sqrt{6}/2} \times S^1_{r=1/\sqrt{2}}$\, is isometric to \,$\R^2/\wt{\Gamma}$\,, and therefore the maximal tori 
\,$(S^1_{r=\sqrt{6}/2} \times S^1_{r=1/\sqrt{2}})/\{\pm \id\}$\, of \,$\SU(3)/\SO(3)$\, are isometric to \,$\R^2/\Gamma$\,, where \,$\Gamma \subset \R^2$\,
is the lattice generated by \,$\wt{\Gamma}$\, and the point \,$(\tfrac{\sqrt{6}}{2}\,\pi, \tfrac{\sqrt{2}}{2}\,\pi)$\, which corresponds to the antipodal
point of the origin in \,$\R^2/\wt{\Gamma}$\,. It can be shown that \,$\Gamma = \Z\,(\tfrac{\sqrt{6}}{2}\,\pi, \tfrac{\sqrt{2}}{2}\,\pi) \oplus
\Z\,(\tfrac{\sqrt{6}}{2}\,\pi, -\tfrac{\sqrt{2}}{2}\,\pi)$\, holds. The two generators of \,$\Gamma$\, 
are not orthogonal to each other (they are at an angle of \,$\tfrac\pi3$\,). It follows that
\,$\R^2/\Gamma$\,, and hence any maximal flat torus of \,$\SU(3)/\SO(3)$\,, is diffeomorphic to \,$S^1 \times S^1$\,, but is not globally
isometric to a product of circles.

\textbf{Type (G).}
The totally geodesic submanifolds corresponding to the Lie triple systems
of type (G) are of course the traces of the geodesics of \,$\SU(3)/\SO(3)$\,; each of them runs within a maximal torus, and their behavior (either they
are periodic, or they are injective and then their trace is dense in the torus) depends on their starting angle in the well-known way.

\textbf{Type (M).}
To construct totally geodesic embeddings of type (M), 
we consider the 3-dimensional complex space \,$V$\, of symmetric complex \,$(2\times 2)$-matrices. \,$V$\, becomes a unitary space via the usual
endomorphism inner product, and this inner product gives rise to the 
Lie group \,$\SU(V) \cong \SU(3)$\,. \,$V$\, has a canonical real form \,$V_{\R} := \Menge{X \in V}{\overline{X}=X}$\,, which we use to define the Lie subgroup
\,$\SO(V) := \Menge{B\in\SU(V)}{B(V_{\R}) = V_{\R}}$\, of \,$\SU(V)$\,,  isomorphic to \,$\SO(3)$\,. 
Thereby we have the realization \,$\SU(V)/\SO(V)$\, of the Riemannian symmetric space \,$\SU(3)/\SO(3)$\,,
and we will construct the totally geodesic submanifolds of type (M) in this realization.

For this, consider the Lie group homomorphism
$$ \Psi: \SU(2) \to \SU(V), \; B \mapsto (X \mapsto BXB^{T}) \;, $$
where \,$B^T$\, denotes the transpose of \,$B$\, for any \,$B \in \SU(2)$\,. Because \,$\SU(2)$\, and \,$\SU(V)$\, are simple Lie groups, \,$\Psi$\, is a homothety
with regard to the invariant Riemannian metrics induced on \,$\SU(2)$\, resp.~\,$\SU(V)$\, by the endomorphism inner product, i.e.~there exists \,$c \in \R_+$\,
so that the linearization \,$\Psi_L: \liesu(2) \to \liesu(V)$\, of \,$\Psi$\, satisfies \,$\|\Psi_L(H)\|^2 = c^2 \,\|H\|^2$\, for all \,$H \in \liesu(2)$\,. 
To determine the value of \,$c$\, we note that \,$\Psi_L$\, is given explicitly by
$$ \Psi_L: \liesu(2) \to \liesu(V),\; H \mapsto (X \mapsto HX + XH^T) \; ; $$
by explicit calculations for \,$H := \left( \begin{smallmatrix} i & 0 \\ 0 & -i \end{smallmatrix} \right) \in \liesu(2)$\, we find
\,$\|H\|^2 = 2$\, and (using the mentioned description of \,$\Psi_L$\,) \,$\|\Psi_L(H)\|^2 = 8$\,. Therefore we have \,$c=2$\,. 

Moreover, we have \,$\Psi^{-1}(\SO(V)) = \SO(2) \cup (J\cdot \SO(2)) =: K$\,, where we again put \,$J := \left( \begin{smallmatrix} i & 0 \\ 0 & -i \end{smallmatrix}
\right) \in \SU(2)$\,, and \,$\Psi$\, is compatible with the involutions on \,$\SU(2)$\, resp.~\,$\SU(V)$\, given by the symmetric space structures of
\,$\SU(2)/\SO(2)$\, resp.~\,$\SU(V)/\SO(V)$\,. Therefore \,$\Psi$\, gives rise to a totally geodesic embedding, which is homothetic with \,$c=2$\,:
$$ \underline{\Psi} : \SU(2)/K \to \SU(V)/\SO(V) \;. $$
The totally geodesic submanifold of \,$\SU(V)/\SO(V)$\, that is the image of \,$\underline{\Psi}$\, turns out to be of type (M). 

As we saw above, \,$\SU(2)/\SO(2)$\, is isometric to \,$S^2_{r=1/\sqrt{2}}$\,, and the pair of points \,$\SO(2)$\,, \,\hbox{$J\cdot \SO(2)$}\, of \,$\SU(2)/\SO(2)$\,
corresponds to a pair of antipodal points in \,$S^2_{r=1/\sqrt{2}}$\, under that isometry. Therefore \,$\SU(2)/K$\, is isometric to \,$\RP^2_{\vkap=2}$\,. 
Because \,$\underline{\Psi}$\, is a homothety with \,$c=2$\,, its image, and hence any totally geodesic submanifold of \,$\SU(3)/\SO(3)$\, of type (M),
is isometric to \,$\RP^2_{\vkap=1/2}$\,.

\end{document}